\documentclass[aop,MSNbibl,secthm,seceqn,dvips]{arximspdf}
\usepackage{mathbh}
\usepackage{graphicx}

\doi{10.1214/10-AOP535}
\volume{38}
\issue{5}
\pubyear{2010}
\firstpage{2066}
\lastpage{2102}

\makeatletter

\newtheorem{lemma}[thm]{Lemma}
\newtheorem{proposition}[thm]{Proposition}
\newtheorem{theorem}[thm]{Theorem}
\newtheorem{conjecture}[thm]{Conjecture}

\makeatother

\begin{document}
\begin{frontmatter}

\title{Intermittency on catalysts: Voter model\protect\thanksref{TIT1}}
\runtitle{Intermittency on catalysts}
\thankstext{TIT1}{Supported in part by
the DFG-NWO Bilateral Research Group ``Mathematical Models from Physics and
Biology'' and by the DFG Research Group 718 ``Analysis and Stochastics
in Complex
Physical Systems.''}

\begin{aug}
\author[A]{\fnms{J.} \snm{G\"artner}\ead[label=e1]{jg@math.tu-berlin.de}},
\author[B]{\fnms{F.} \snm{den Hollander}\ead[label=e2]{denholla@math.leidenuniv.nl}\corref{}}
and
\author[C]{\fnms{G.} \snm{Maillard}\thanksref{tit2}\ead[label=e3]{maillard@cmi.univ-mrs.fr}}
\runauthor{J. G\"artner, F. den Hollander and G. Maillard}
\thankstext{tit2}{Supported by a postdoctoral fellowship from The
Netherlands Organization for Scientific Research (Grant 613.000.307) while
at EURANDOM.}
\affiliation{Technische Universit\"at Berlin, Leiden University and
EURANDOM and Universit\'e de Provence}
\address[A]{J. G\"artner\\
Institut f\"ur Mathematik\\
Technische Universit\"at Berlin\\
Strasse des 17. Juni 136\\
D-10623 Berlin\\
Germany\\
\printead{e1}} 
\address[B]{F. den Hollander\\
Mathematical Institute\\
Leiden University\\
P.O. Box 9512\\
2300 RA Leiden\\
The Netherlands\\
and\\
EURANDOM\\
P.O. Box 513\\
5600 MB Eindhoven\\
The Netherlands\\
\printead{e2}}
\address[C]{G. Maillard\\
CMI-LATP\\
Universit\'e de Provence\\
39 rue F. Joliot-Curie\\
F-13453 Marseille Cedex 13\\
France\\
\printead{e3}}
\end{aug}

\received{\smonth{8} \syear{2009}}
\revised{\smonth{1} \syear{2010}}

%
\begin{abstract}
In this paper we study intermittency for the parabolic Anderson equation
$\partial u/\partial t = \kappa\Delta u + \gamma\xi u$ with $u\dvtx
\mathbb{Z}^d
\times[0,\infty)\to\mathbb{R}$, where $\kappa\in[0,\infty)$ is
the diffusion constant,
$\Delta$ is the discrete Laplacian, $\gamma\in(0,\infty)$ is the
coupling constant,
and $\xi\dvtx\mathbb{Z}^d\times[0,\infty)\to\mathbb{R}$ is a
space--time random medium. The
solution of this equation describes the evolution of a ``reactant'' $u$
under the
influence of a ``catalyst'' $\xi$.

We focus on the case where $\xi$ is the voter model with
opinions 0 and 1 that are updated according to a random walk transition kernel,
starting from either the Bernoulli measure $\nu_\rho$ or the
equilibrium measure
$\mu_\rho$, where $\rho\in(0,1)$ is the density of 1's. We consider
the annealed
Lyapunov exponents, that is, the exponential growth rates of the
successive moments
of $u$. We show that if the random walk transition kernel has zero mean
and finite
variance, then these exponents are trivial for $1\leq d\leq4$, but
display an
interesting dependence on the diffusion constant $\kappa$ for $d\geq
5$, with
qualitatively different behavior in different dimensions.

In earlier work we considered the case where $\xi$ is a field of independent
simple random walks in a Poisson equilibrium, respectively, a symmetric
exclusion
process in a Bernoulli equilibrium, which are both reversible dynamics.
In the
present work a main obstacle is the nonreversibility of the voter model
dynamics,
since this precludes the application of spectral techniques. The
duality with
coalescing random walks is key to our analysis, and leads to a representation
formula for the Lyapunov exponents that allows for the application of large
deviation estimates.
\end{abstract}

%

\end{frontmatter}

\section{Introduction and main results}
\label{S1}

The outline of this section is as follows. In Section \ref{S1.1} we provide
motivation. In Sections \ref{S1.2}--\ref{S1.4} we recall some basic facts
about the voter model. In Section \ref{S1.5} we define the annealed Lyapunov
exponents, which are the main objects of our study. In Section~\ref{S1.6}
we prove a representation formula for these exponents in terms of coalescing
random walks released at Poisson times along a random walk path. This
representation formula is the starting point for our further analysis.
Our main theorems are stated in Section \ref{S1.7} (and proved in
Sections \ref{S2}--\ref{S5}). Finally, in Sections \ref{S1.8}--\ref{S1.9}
we list some open problems and state a scaling conjecture.

\subsection{Reactant and catalyst}
\label{S1.1}

The \textit{parabolic Anderson equation} is the partial differential equation
%
\begin{equation}\label{pA}
\frac{\partial}{\partial t}u(x,t) = \kappa\Delta u(x,t) + \gamma\xi
(x,t)u(x,t),\qquad
x\in\mathbb{Z}^d, t\geq0.
\end{equation}
Here, the $u$-field is $\mathbb{R}$-valued, $\kappa\in[0,\infty)$
is the diffusion
constant, $\Delta$ is the discrete Laplacian, acting on $u$ as
%
\begin{equation}\label{dL}
\Delta u(x,t) = \mathop{\sum_{y\in\mathbb{Z}^d }}_{\|y-x\|=1} [u(y,t)-u(x,t)]
\end{equation}
($\|\cdot\|$ is the Euclidean norm), $\gamma\in[0,\infty)$ is the
coupling constant,
while
%
\begin{equation}\label{rf}
\xi= \{\xi(x,t) \dvtx x\in\mathbb{Z}^d, t\geq0\}
\end{equation}
is an $\mathbb{R}$-valued random field that evolves with time and that drives
the equation.
As initial condition for (\ref{pA}) we take
%
\begin{equation}\label{ic}
u(\cdot,0) \equiv1.
\end{equation}

The PDE in (\ref{pA}) describes the evolution of a system of two types
of particles,
$A$~and $B$, where the $A$-particles perform autonomous dynamics and
the $B$-particles
perform independent simple random walks that branch at a rate that is
equal to $\gamma$
times the number of $A$-particles present at the same location. The
link is that $u(x,t)$
equals the average number of $B$-particles at site $x$ at time $t$
conditioned on the
evolution of the $A$-particles. The initial condition in (\ref{ic})
corresponds to
starting off with one $B$-particle at each site. Thus, the solution of
(\ref{pA}) may
be viewed as describing the evolution of a \textit{reactant} $u$ under
the influence of a
\textit{catalyst} $\xi$. Our focus of interest will be on the \emph
{annealed Lyapunov
exponents}, that is, the exponential growth rates of the successive
moments of $u$.

In earlier work (G\"artner and den Hollander \cite{garhol06}, G\"artner, den
Hollander and Maillard \cite{garholmai07,garholmai08pr}) we
treated the
case where $\xi$ is a field of independent simple random walks in a
Poisson equilibrium,
respectively, a symmetric exclusion process in a Bernoulli equilibrium.
In the present
paper we focus on the case where $\xi$ is the \textit{Voter Model}
(VM), that is, $\xi$ takes
values in $\{0,1\}^{\mathbb{Z}^d \times[0,\infty)}$, where $\xi
(x,t)$ is
the opinion of
site $x$ at time $t$, and opinions are imposed according to a random
walk transition
kernel. We choose $\xi(\cdot,0)$ according to either the Bernoulli
measure $\nu_\rho$
or the equilibrium measure $\mu_\rho$, where $\rho\in(0,1)$ is the
density of $1$'s.
We may think of 0 as a vacancy and 1 as a particle.

An overview of the main results in \cite{garhol06,garholmai07,garholmai08pr}
and the present paper as well as further literature is given in G\"
artner, den Hollander
and Maillard \cite{garholmai08HvW}.
G\"artner and Heydenreich \cite{garhey06} consider the case where the
catalyst consists of
a single random walk.

\subsection{Voter model}\label{S1.2}

Throughout the paper we abbreviate $\Omega=\{0,1\}^{\mathbb{Z}^d}$\break (equipped
with the product
topology), and we let $p\dvtx\mathbb{Z}^d \times\mathbb{Z}^d\to
[0,1]$ be the
transition kernel of an
irreducible random walk, that is,
%
\begin{eqnarray}\label{pdef}
&\displaystyle\sum_{y\in\mathbb{Z}^d} p(x,y)=1 \qquad \forall x\in\mathbb{Z}^d,\nonumber\\
&\hspace*{96pt}p(x,y)=p(0,y-x) \geq0 \qquad \forall x,y\in\mathbb{Z}^d,\\
&p(\cdot,\cdot) \mbox{ generates } \mathbb{Z}^d.\nonumber
\end{eqnarray}
Occasionally we will need to assume that $p(\cdot,\cdot)$ has zero
mean and finite variance.
A special case is simple random walk
%
\begin{equation}\label{SRW}
p(x,y) =
\cases{
\displaystyle\frac{1}{2d}, &\quad  if $\|x-y\|=1$,\vspace*{2pt}\cr
0, &\quad  otherwise.
}
\end{equation}

The VM is the Markov process on $\Omega$ whose generator $L$ acts
on cylindrical functions $f$ as
%
\begin{equation}\label{VMgen}
(Lf)(\eta) = \sum_{x,y\in\mathbb{Z}^d}
p(x,y) [f(\eta^{x\rightarrow y})- f(\eta) ],\qquad  \eta\in\Omega,
\end{equation}
where
%
\begin{equation}\label{VMtr}
\eta^{x\rightarrow y}(z) =
\cases{
\eta(x), &\quad if $ z=y$,\cr
\eta(z), &\quad if  $z\neq y$.
}
\end{equation}
Under this dynamics, site $x$ imposes its state on site $y$ at rate
$p(x,y)$. The states
$0$ and $1$ are referred to as opinions or, alternatively, as vacancy
and particle. The
VM is a \textit{nonconservative} dynamics: opinions are not preserved.
We write
$(S_t)_{t\geq0}$ to denote the Markov semigroup associated with $L$.

Let $\xi_t = \{\xi(x,t); x\in\mathbb{Z}^d\}$ be the random
configuration of
the VM at time
$t$. Let $\mathbb{P}_\eta$ denote the law of $\xi$ starting from
$\xi
_0=\eta$, and let $\mathbb{P}_\mu
= \int_\Omega\mu(d\eta) \mathbb{P}_\eta$. We will consider two choices
for the starting measure
$\mu$:
%
\begin{equation}\label{muchoice}
\cases{
\mu= \nu_\rho,&\quad  the Bernoulli measure with density $\rho\in
(0,1)$,\cr
\mu= \mu_\rho,&\quad  the equilibrium measure with density $\rho
\in(0,1)$.
}
\end{equation}

Let $p^*(\cdot,\cdot)$ be the \textit{dual} transition kernel, defined
by $p^*(x,y)=p(y,x)$,
$x,y\in\mathbb{Z}^d$, and $p^{(s)}(\cdot,\cdot)$ the \textit{symmetrized}
transition kernel, defined
by $p^{(s)}(x,y) = (1/2)[p(x,y)+p^*(x,y)]$, $x,y\in\mathbb{Z}^d$.
The ergodic properties of the
VM are qualitatively different for recurrent and for transient
$p^{(s)}(\cdot,\cdot)$.
In particular, when $p^{(s)}(\cdot,\cdot)$ is
\textit{recurrent} all equilibria are \textit{trivial}, that is, $\mu
_\rho= (1-\rho)\delta_0
+\rho\delta_1$, while when $p^{(s)}(\cdot,\cdot)$ is \emph
{transient} there are also
\textit{nontrivial} equilibria, that is, ergodic measures $\mu_\rho$.
In the latter case,
$\mu_\rho$ is taken to be the unique shift-invariant and ergodic
equilibrium with density
$\rho$. For both cases we have
%
\begin{equation}\label{VMerg1}
\mathbb{P}_{\nu_\rho} (\xi_t\in\cdot) \to\mu_\rho(\cdot)
\qquad \mbox{weakly as } t\to\infty,
\end{equation}
with the same convergence for any starting measure $\mu$ that is
stationary and
ergodic with density $\rho$ (see Liggett \cite{lig85}, Corollary V.1.13).

\begin{figure}[b]

\includegraphics{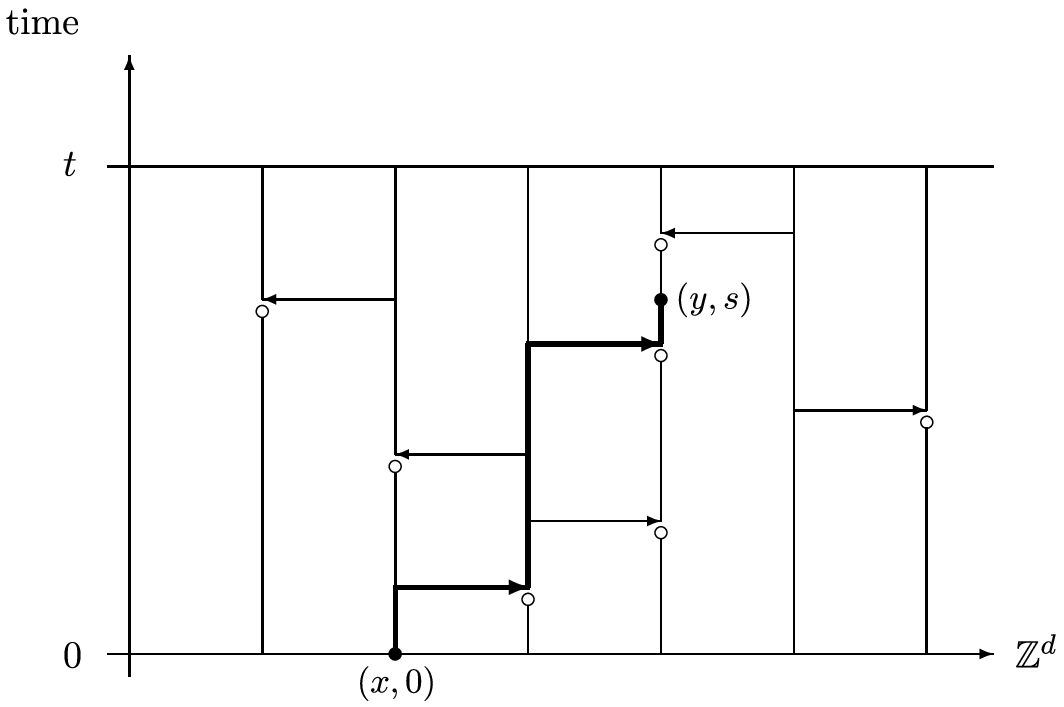}

  \caption{Graphical representation $\mathcal{G}_t$. Opinions propagate
along paths.}\label{fig1}
\end{figure}

We will frequently use the measures $\nu_\rho S_T$, $T\in[0,\infty
]$, where $\nu_\rho
S_\infty=\mu_\rho$ by convention in view of (\ref{VMerg1}). The VM
is \textit{attractive}
(see Liggett \cite{lig85}, Definition~III.2.1 and Theorem III.2.2).
Consequently, since
$\nu_\rho$ has positive correlations, the same is true for $\nu_\rho
S_T$, that is,
nondecreasing functions on $\Omega$ are positively correlated (see
Liggett \cite{lig85},
Theorem II.2.14).

\subsection{Graphical representation and duality}
\label{S1.3}

In the VM's \textit{graphical representation} $\mathcal{G}_t$ from time
$0$ up
to time $t$
(see, e.g., Cox and Griffeath \cite{coxgri83}, Section 0), space is
drawn sideward,
time is drawn upward, and for each ordered pair of sites $x,y\in
\mathbb{Z}^d$
arrows are
drawn from $x$ to $y$ at Poisson rate $p(x,y)$. A path from $(x,0)$ to $(y,s)$,
$s\in(0,t]$, in $\mathcal{G}_t$ (see Figure \ref{fig1}) is a sequence of
space--time points
$(x_0,s_0),(x_0,s_1),(x_1,s_1),\dots,$ $(x_n,s_n),(x_n,s_{n+1})$ such that:
\begin{longlist}[(iii)]
\item[(i)] $x_0=x,s_0=0,x_n=y,s_{n+1}=s$;
\item[(ii)] the sequence of times $(s_i)_{0\leq i \leq n+1}$ is increasing;
\item[(iii)] for each $1\leq i\leq n$, there is an arrow from
$(x_{i-1},s_{i})$ to
$(x_i,s_i)$;
\item[(iv)] for each $0\leq i\leq n$, no arrow \textit{points to} $x_i$
at any time
in $(s_i,s_{i+1})$.
\end{longlist}
Then $\xi$ can be represented as
%
\begin{equation}\label{vmgraph}
\qquad \xi(y,s)=
\cases{
1, &\quad if there exists a path from $(x,0)$ to $(y,s)$ in $\mathcal
{G}_t$\cr
&\quad for some $x\in\xi(0)$,\cr
0, &\quad otherwise,
}
\end{equation}
where $\xi(0)=\{x\in\mathbb{Z}^d\dvtx\xi(x,0)=1\}$ is the set of initial
locations of the
1's. The graphical representation corresponds to binary \textit{branching
} with
transition kernel $p(\cdot,\cdot)$ and step rate $1$ and \textit{killing
} at the
moment when an arrow comes in from another location. Figure \ref{fig1} shows how
opinions propagate along paths. An open circle indicates that the site
adopts the opinion of the site where the incoming arrow comes from.
The thick line from $(x,0)$ to $(y,s)$ shows that the opinion
at site $y$ at time $s$ stems from the opinion at a unique site $x$ at
time $0$.

%
%
%
%
%
%
%
%

\begin{figure}[b]

\includegraphics{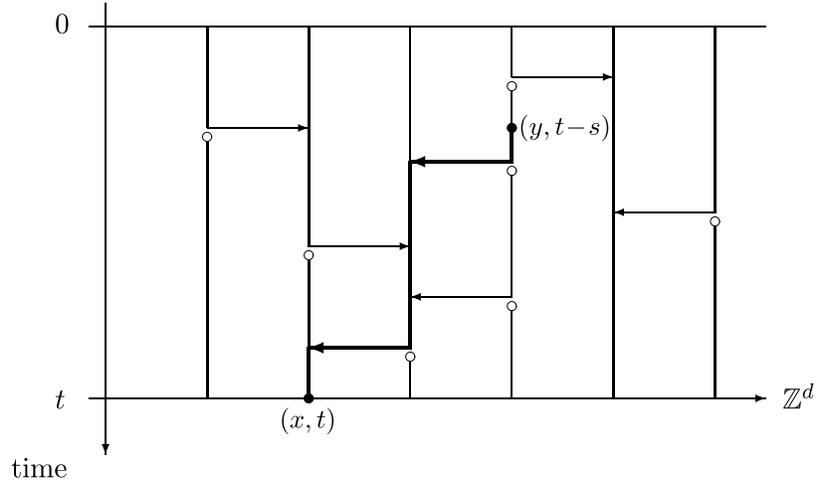}

  \caption{Dual graphical representation $\mathcal{G}_t^\ast$.
Opinions propagate along time-reversed coalescing paths.}\label{fig2}
\end{figure}

We can define the \textit{dual graphical representation} ${\mathcal
{G}}_t^\ast$
by reversing time and direction of all the arrows in ${\mathcal{G}}_t$.
The dual process $(\xi^\ast_s)_{0\leq s \leq t}$ on $\mathcal
{G}^\ast_t$
can then
be represented as
%
\begin{equation}\label{coalgraph}
\qquad \xi^\ast(x,t)=
\cases{
1, &\quad  if there exists a path from $(y,t-s)$ to $(x,t)$ in $\mathcal{G}
^\ast_t$\cr
&\quad  for some $y\in\xi^\ast(t-s)$,\cr
0, &\quad otherwise,
}
\end{equation}
where $\xi^\ast(t-s)=\{x\in\mathbb{Z}^d\dvtx \xi^\ast(x,t-s)=1\}
$. The
dual graphical representation
corresponds to \textit{coalescing random walks} with dual transition
kernel $p^*(\cdot,\cdot)$
and step rate $1$ (see Figure \ref{fig2}).


%
%
%
%
%
%
%
%
%
%


Figures \ref{fig1} and \ref{fig2} make it plausible that the equilibrium measure $\mu
_\rho$ in (\ref{VMerg1})
is \textit{nonreversible}, because the evolution is not invariant under
time reversal.


\subsection{Correlation functions}
\label{S1.4}

A key tool in the present paper is the following \textit{representation
formula for the
$n$-point correlation functions} of the VM, which is an immediate
consequence of the
dual graphical re\-presentation (see, e.g., Cox and Griffeath \cite
{coxgri83}, Section 1).
For $n \in\mathbb{N}$, $x_1,\dots,x_n\in\mathbb{Z}^d$ and $-\infty
< s_1 \leq
\cdots\leq s_n \leq t$, let
%
\begin{equation}\label{coaldef}
\xi^\ast_t \{(x_1,s_1),\dots,(x_n,s_n) \}
\end{equation}
be the set of locations at time $t$ of $n$ coalescing random walks,
with transition
kernel $p^*(\cdot,\cdot)$ and step rate 1, when the $m$th random walk
is born at site
$x_m$ at time $s_m$, $1 \leq m \leq n$, and let
%
\begin{equation}\label{Ntdef}
\mathcal{N}_t \{(x_1,s_1),\dots,(x_n,s_n) \}
= | \xi^\ast_t \{(x_1,s_1),\dots,(x_n,s_n) \} |\\
\end{equation}
be the number of random walks alive at time $t$.

The following lemma gives us a handle on the $n$-point correlation functions.

\begin{lemma}\label{reprform}
For all $n\in\mathbb{N}$, $T\in[0,\infty]$, $x_1,\dots,x_n\in
\mathbb{Z}^d$
and $-\infty<s_1\leq\cdots\leq s_n\leq t<\infty$,
%
\begin{equation}\label{dual-vm}
\qquad \mathbb{P}_{\nu_\rho S_T} \bigl(\xi(x_m,t-s_m)=1 \ \forall1\leq m\leq n \bigr)
= \mathbb{E}^\ast\bigl(\rho^{ \mathcal{N}_{T+t} \{(x_1,s_1),\dots
,(x_n,s_n) \}} \bigr),
\end{equation}
where $\mathbb{E}^\ast$ denotes expectation with respect to the coalescing
random walk
dynamics.
\end{lemma}

\begin{pf}
For $T<\infty$, we have
\begin{eqnarray}\label{dual-vm0}
&&
\mathbb{P}_{\nu_\rho S_T} \bigl(\xi(x_m,t-s_m)=1 \ \forall1\leq m\leq n
\bigr)\nonumber\\[-8pt]\\[-8pt]
&&\qquad  =\mathbb{P}_{\nu_\rho} \bigl(\xi(x_m,T+t-s_m)=1 \ \forall1\leq m\leq n \bigr).\nonumber
\end{eqnarray}
The event in the right-hand side of (\ref{dual-vm0}) occurs if and
only if $\xi(z,0)=1$
for all sites $z$ in the set $\xi^\ast_{T+t}\{(x_1,s_1),\dots
,(x_n,s_n)\}$ (Figure \ref{fig2}),
which under $\nu_\rho$ has probability $\rho^{ \mathcal{N}_{T+t}\{
(x_1,s_1),\dots,(x_n,s_n)\}}$
and proves the claim. Since $t\mapsto\mathcal{N}_t$ is nonincreasing,
we may
let $T\to\infty$
in (\ref{dual-vm}) and use (\ref{VMerg1}) to get the formula for
$T=\infty$.
\end{pf}

Note that for $T=\infty$ the right-hand side of (\ref{dual-vm}) does
not depend on
$t$, in accordance with the fact that $\nu_\rho S_\infty= \mu_\rho$
is an
equilibrium measure.


\subsection{Lyapunov exponents}
\label{S1.5}

By the Feynman--Kac formula, the formal solution of (\ref{pA}) and
(\ref{ic})
reads
%
\begin{equation}\label{fey-kac1}
u(x,t) = \mathrm{E}_{ x} \biggl(\exp\biggl[\gamma\int_0^t
\xi\bigl(X^\kappa(s),t-s \bigr)\,ds \biggr] \biggr),
\end{equation}
where $X^\kappa$ is a simple random walk on $\mathbb{Z}^d$ with step rate
$2d\kappa$,
and $\mathrm{E}_{ x}$ denotes expectation w.r.t. $X^\kappa$ given
$X^\kappa(0)=x$.
Let $\mu$ be an arbitrary initial distribution.
For $p\in\mathbb{N}$ and $t>0$, the $p$th moment of the solution is then
given by
%
\begin{equation}\label{fkmom}
\qquad \mathbb{E}_{ \mu}([u(0,t)]^p) = (\mathbb{E}_{ \mu} \otimes\mathrm
{E}_{ 0}^{\otimes p})
\biggl(\exp\biggl[\gamma\int_0^t \sum_{q=1}^{p}
\xi\bigl(X_q^\kappa(s),t-s \bigr) \,ds \biggr] \biggr),
\end{equation}
where $X_q^\kappa$, $q=1,\dots,p$, are $p$ independent copies of
$X^\kappa$.

For $p\in\mathbb{N}$ and $t>0$, define
%
\begin{equation}\label{lyapdef}
\Lambda^\mu_p(t) = \frac{1}{pt} \log\mathbb{E}_{ \mu}([u(0,t)]^p).
\end{equation}
Then
%
\begin{equation}\label{fey-kac2}
\Lambda^\mu_p(t) = \frac{1}{pt} \log(\mathbb{E}_{ \mu} \otimes
\mathrm{E}_{
0}^{\otimes p})
\Biggl(\exp\Biggl[\gamma\int_0^t \sum_{q=1}^{p}
\xi\bigl(X_q^\kappa(s),t-s \bigr)\,ds \Biggr] \Biggr).
\end{equation}
We will see that for $\mu=\nu_\rho S_T$, $T\in[0,\infty]$, the last quantity
admits a limit as $t\to\infty$,
%
\begin{equation}\label{lyap2}
\lambda^\mu_p = \lim_{t\to\infty} \Lambda^\mu_p(t),
\end{equation}
which is independent of $T$ and which we call the $p$th \textit{annealed
Lyapunov exponent}.
Note that $\Lambda^\mu_p(t) \in[\rho\gamma,\gamma]$ for all
$t>0$, as is immediate from
(\ref{fey-kac2}) and Jensen's inequality. Hence,
%
\begin{equation}\label{lyasand}
\lambda^\mu_p \in[\rho\gamma,\gamma].
\end{equation}

From H\"older's inequality applied to (\ref{lyapdef}), it follows that
$\Lambda^\mu_p(t)
\geq\Lambda^\mu_{p-1}(t)$ for all $t>0$ and $p\in\mathbb
{N}\setminus\{1\}
$. Hence, $\lambda^\mu_p
\geq\lambda^\mu_{p-1}$ for all $p\in\mathbb{N}\setminus\{1\}$.
We say
that the solution of the
parabolic Anderson model is \textit{$p$-intermittent} if $\lambda^\mu
_p>\lambda^\mu_{p-1}$.
In the latter case the solution is $q$-intermittent for all $q>p$ as
well (see, e.g., G\"artner
and Heydenreich \cite{garhey06}, Lemma 3.1). We say that the solution
is \textit{intermittent}
if it is $p$-intermittent for all $p\in\mathbb{N}\setminus\{1\}$.
Intermittent means that the
$u$-field develops sparse high peaks dominating the moments in such a
way that each moment
is dominated by its own collection of peaks (see G\"artner and K\"onig
\cite{garkon05},
Section 1.3, and G\"artner and den Hollander \cite{garhol06}, Section 1.2).

\subsection{Representation formula}
\label{S1.6}

In this section we derive a \textit{coalescing random walk
representation for
the Lyapunov exponents}. Recall (\ref{Ntdef}). For $n \in\mathbb{N}$,
$x_1,\dots,x_n\in\mathbb{Z}^d$ and $-\infty< s_1 \leq\cdots\leq
s_n \leq
t$, let
%
\begin{equation}\label{Ntcoaldef}
\mathcal{N}_t^{\mathrm{coal}}\{(x_1,s_1),\dots,(x_n,s_n) \}
= n - \mathcal{N}_t \{(x_1,s_1),\dots,(x_n,s_n) \}
\end{equation}
be the number of random walks coalesced at time $t$.
Let $\Pi_{\rho\gamma}$ and $\mathbb{P}_{\mathrm{Poiss}}$ denote the
Poisson point
process on
$\mathbb{R}$ with intensity $\rho\gamma$ and its law, respectively.
We consider $\Pi_{\rho\gamma}$ as a random subset of $\mathbb{R}$
and write
$\Pi_{\rho\gamma}(B)=\Pi_{\rho\gamma}\cap B$ for Borel sets
$B\subseteq\mathbb{R}$.

\begin{proposition}\label{MomCRWprop}
For all $T\in[0,\infty]$, $t>0$ and right-continuous paths
$\varphi_q\dvtx[0,t]\to\mathbb{Z}^d$, $q=1,\dots,p$,
\begin{eqnarray}\label{momCRW}
&&e^{-\rho\gamma pt} \mathbb{E}_{ \nu_\rho S_T}
\Biggl(\exp\Biggl[\gamma\int_0^t \sum_{q=1}^{p}
\xi\bigl(\varphi_q(s),t-s \bigr) \,ds \Biggr]
\Biggr)\nonumber\\[-8pt]\\[-8pt]
&&\qquad  = (\mathbb{E}_{\mathrm{Poiss}}^{\otimes p}\otimes\mathbb{E}^\ast)
\bigl(\rho^{-\mathcal{N}^{\mathrm{coal}}_{T+t} \{\bigcup_{q=1}^{p} \{
(\varphi_q(s),s)
\dvtx s\in\Pi_{\rho\gamma}^{(q)}([0,t]) \} \}} \bigr),\nonumber
\end{eqnarray}
where $\Pi_{\rho\gamma}^{(q)}$, $q=1,\dots,p$, are $p$ independent
copies of
$\Pi_{\rho\gamma}$. In particular,
\begin{eqnarray}\label{lyaCRW}
&&\exp\bigl[pt\bigl(\Lambda^{\nu_\rho S_T}_p(t)-\rho\gamma\bigr)
\bigr]\nonumber\\[-8pt]\\[-8pt]
&&\qquad  = (\mathrm{E}_{ 0}^{\otimes p}\otimes\mathbb{E}_{\mathrm
{Poiss}}^{\otimes p}\otimes\mathbb{E}^\ast)
\bigl(\rho^{-\mathcal{N}^{\mathrm{coal}}_{T+t} \{\bigcup_{q=1}^{p} \{
(X_q^\kappa(s),s)
\dvtx s\in\Pi_{\rho\gamma}^{(q)}([0,t]) \} \}} \bigr).\nonumber
\end{eqnarray}
\end{proposition}

\begin{pf}
Fix $\varphi_q$, $q=1,\dots,p$. By a Taylor expansion of the factors
 $\exp[\gamma\times\int_0^t \xi(\varphi_q(s),t-s) \,ds]$,
$q=1,\dots,p$,
we have
%
\begin{eqnarray}\label{lyaCRW3}
&&e^{-\rho\gamma pt} \mathbb{E}_{ \nu_\rho S_T}
\Biggl(\exp\Biggl[\gamma\int_0^t \sum_{q=1}^{p}
\xi\bigl(\varphi_q(s),t-s \bigr) \,ds \Biggr] \Biggr)\nonumber\\
&&\qquad  =e^{-\rho\gamma pt} \Biggl[\prod_{q=1}^p\sum_{n_q=0}^{\infty}
\frac{\gamma^{n_q}}{n_q!}
\Biggl(\prod_{m=1}^{n_q} \int_0^t ds_m^{(q)} \Biggr) \Biggr] \nonumber\\
&&\qquad \quad {}
\times\mathbb{E}_{ \nu_\rho S_T} \Biggl(\prod_{q=1}^{p}\prod_{m=1}^{n_q}
\xi\bigl(\varphi_q \bigl(s_m^{(q)} \bigr),t-s_m^{(q)} \bigr) \Biggr)\\
&&\qquad  = \Biggl[\prod_{q=1}^p\sum_{n_q=0}^{\infty}
\frac{(\rho\gamma t)^{n_q}}{n_q!} e^{-\rho\gamma t}
\frac{1}{t^{n_q}} \Biggl(\prod_{m=1}^{n_q}\int_0^t ds_m^{(q)} \Biggr) \Biggr]\nonumber\\
&&\qquad \quad {}
\times\rho^{-\sum_{q=1}^p n_q}
\mathbb{E}_{ \nu_\rho S_T} \Biggl(\prod_{q=1}^{p}\prod_{m=1}^{n_q}
\xi\bigl(\varphi_q \bigl(s_m^{(q)} \bigr),t-s_m^{(q)} \bigr) \Biggr).\nonumber
\end{eqnarray}
For each $q=1,\ldots,p$:
\begin{itemize}
\item
$[(\rho\gamma t)^{n_q}/n_q!]\exp[-\rho\gamma t]$, $n_q\in\mathbb
{N}_0=\mathbb{N}
\cup\{0\}$,
is the Poisson distribution with parameter $\rho\gamma t$;
\item
$(1/t^{n_q})(\prod_{m=1}^{n_q} \int_0^t ds_m^{(q)})$ is the uniform
distribution
on $[0,t]^{n_q}$,
coinciding with the distribution of the (unordered) points of
$\Pi^{(q)}_{\rho\gamma}$ in $[0,t]$ given
$ |\Pi^{(q)}_{\rho\gamma}([0,t]) |=n_q$, $n_q\in\mathbb{N}_0$.
\end{itemize}
Moreover, by Lemma \ref{reprform}, we have
\begin{eqnarray}\label{lyaCRW5}
&&\mathbb{E}_{ \nu_\rho S_T} \Biggl(\prod_{q=1}^{p}\prod_{m=1}^{n_q}
\xi\bigl(\varphi_q \bigl(s_m^{(q)} \bigr),t-s_m^{(q)} \bigr) \Biggr)\nonumber\\[-8pt]\\[-8pt]
&& \qquad  =\mathbb{E}^\ast\bigl(\rho^{ \mathcal{N}_{T+t} \{\bigcup_{q=1}^{p} \{
(\varphi_q(s_m^{(q)}),s_m^{(q)} )\dvtx m=1,\dots,n_q \} \}} \bigr).\nonumber
\end{eqnarray}
Therefore, combining (\ref{lyaCRW3}) and (\ref{lyaCRW5}) and inserting
(\ref{Ntcoaldef}),
we get (\ref{momCRW}).

Recalling (\ref{fey-kac2}), we see that formula (\ref{lyaCRW})
follows from
(\ref{momCRW}) by substituting $\varphi_q=X^\kappa_q$, $q=1,\dots
,p$, and
taking the expectation $\mathrm{E}_{ 0}^{\otimes p}$.
\end{pf}

What (\ref{lyaCRW}) in Proposition \ref{MomCRWprop} says is that, for initial
distribution $\mu=\nu_\rho S_T$, the $p$th Lyapunov exponent $\lambda
^\mu_p$
can be computed by taking $p$ simple random walks (with step rate
$2d\kappa$),
releasing coalescing random walks [with dual transition kernel
$p^*(\cdot,\cdot)$ and step rate $1$] from the paths of these $p$
random walks
at rate $\rho\gamma$ until time $t$, recording the total number of
coalescences
up to time $T+t$, and letting $t\to\infty$ afterward. The representation
formula (\ref{lyaCRW}) will be the starting point of our large deviation
analysis.


\subsection{Main theorems}
\label{S1.7}

Theorems \ref{Lyaexist}--\ref{Lyaprop} below are our main results. We write
$\lambda^\mu_p(\kappa)$ to exhibit the $\kappa$-dependence of the Lyapunov
exponents $\lambda^\mu_p$.
The dependence on the other parameters will generally be suppressed
from the
notation.

\begin{theorem}\label{Lyaexist}
For all $d\geq1$, $p\in\mathbb{N}$, $\kappa\in[0,\infty)$,
$\gamma\in
(0,\infty)$ and
$\rho\in(0,1)$, the limit $\lambda^\mu_p$ in (\ref{lyap2})
exists for
$\mu=\nu_\rho S_T$ and is the same for all $T\in[0,\infty]$ (and is
henceforth
denoted by $\lambda_p$).
\end{theorem}

\begin{theorem}\label{Lyalip}
For all $d\geq1$, $p\in\mathbb{N}$, $\gamma\in(0,\infty)$ and
$\rho\in(0,1)$:

\begin{enumerate}[(ii)]
\item[(i)]
$\kappa\mapsto\lambda_p(\kappa)$ is globally Lipschitz outside any
neighborhood
of $0$;
\item[(ii)]
$\lambda_p(\kappa)>\rho\gamma$ for all $\kappa\in[0,\infty)$.
\end{enumerate}
\end{theorem}

\begin{theorem}\label{Lyaprop}
Fix $p\in\mathbb{N}$, $\gamma\in(0,\infty)$ and $\rho\in(0,1)$.

\begin{enumerate}[(ii)]
\item[(i)] If $1\leq d\leq4$ and $p(\cdot,\cdot)$ has zero mean and finite
variance, then $\lambda_p(\kappa)=\gamma$ for all $\kappa\in
[0,\infty)$.
\item[(ii)] If $d\geq5$, then:

\begin{enumerate}[(a)]
\item[(a)] $\lim_{\kappa\downarrow0}\lambda_p(\kappa)=\lambda_p(0)$;
\item[(b)] $\lim_{\kappa\to\infty}\lambda_p(\kappa)=\rho\gamma$;
\item[(c)] if $p(\cdot,\cdot)$ has zero mean and finite variance, then
there exists $\kappa_0>0$ such that $p\mapsto\lambda_p(\kappa)$
is strictly increasing for $\kappa\in[0,\kappa_0)$.
\end{enumerate}
\end{enumerate}
\end{theorem}


%
%

\begin{figure}

\includegraphics{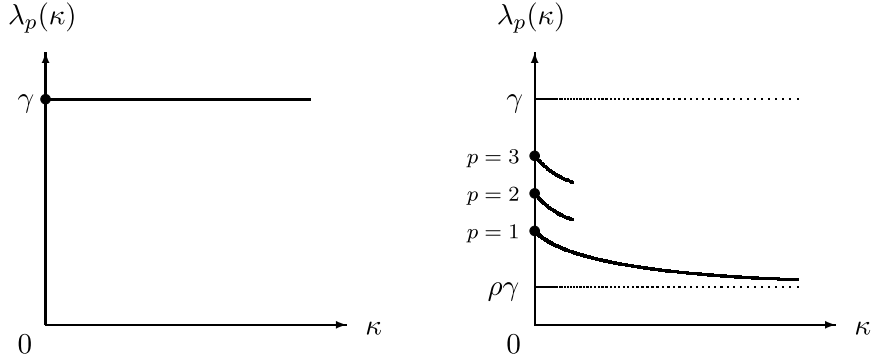}

  \caption{$\kappa\mapsto\lambda_p(\kappa)$ for $1 \leq d \leq4$,
respectively, $d\geq5$, when $p(\cdot,\cdot)$ has zero
mean and finite variance.}\label{fig3}
\end{figure}

Theorem \ref{Lyaexist} says that the Lyapunov exponents exist and do
not depend
on the choice of the starting measure $\mu$. Theorem \ref{Lyalip}
says that the
Lyapunov exponents are continuous functions of the diffusion constant
$\kappa$
away from $0$ and that the system exhibits \textit{clumping} for all
$\kappa$:
the Lyapunov exponents are strictly larger in the random medium than in the
average medium. Theorem \ref{Lyaprop} shows that the Lyapunov exponents
satisfy a \textit{dichotomy} (see Figure \ref{fig3}): for $p(\cdot,\cdot)$ with
zero mean
and finite variance they are trivial when $1 \leq d \leq4$, but
display an
interesting dependence on $\kappa$ when $d\geq5$. In the latter case
(a) the
Lyapunov exponents are continuous in $\kappa$ at $\kappa=0$; (b) the clumping
vanishes in the limit as $\kappa\to\infty$: when the reactant
particles move
much faster than the catalyst particles, they effectively see the
average medium;
(c) the system is intermittent for small $\kappa$: when the reactant particles
move much slower than the catalyst particles, the growth rates of their
successive
moments are determined by different piles of the catalyst.

Theorems \ref{Lyaexist} and \ref{Lyalip} are proved in Sections \ref
{S2} and
\ref{S3}, respectively. Section~\ref{S4} contains block estimates for
coalescing
random walks, which are needed to exploit Proposition \ref{MomCRWprop}
in order
to prove Theorem \ref{Lyaprop}(ii)(a) and (b).
Finally, Theorem~\ref{Lyaprop}(i) and (ii)(c) is proved in Section~\ref{S5}.


\subsection{Open problems}
\label{S1.8}

The following problems remain open:
\begin{enumerate}[(2)]
\item[(1)]
Show that $\lambda_p(\kappa)<\gamma$ for \textit{all } $\kappa\in
[0,\infty)$
when $d\geq5$ and $p(\cdot,\cdot)$ has zero mean and finite variance.
\item[(2)]
Show that $\kappa\mapsto\lambda_p(\kappa)$ is convex on $[0,\infty
)$. Convexity,
when combined with the properties in Theorems \ref{Lyalip}(ii) and
\ref{Lyaprop}(ii)(b),
would imply that $\kappa\mapsto\lambda_p(\kappa)$ is strictly
decreasing on $[0,\infty)$
when $d\geq5$. Convexity was
proved in \cite{garhol06} and \cite{garholmai07} for the case where
$\xi$ is a field
of independent simple random walks in a Poisson equilibrium,
respectively, a symmetric
exclusion process in a Bernoulli equilibrium.
\item[(3)]
Show that the following extension of Theorem \ref{Lyaprop} is true: the
Lyapunov exponents are nontrivial if and only if $p^{(s)}(\cdot,\cdot
)$ is
strongly transient, that is, $\int_0^\infty t p_t^{(s)}(0,0) \,dt
<\infty$. A similar
\textit{full dichotomy} was found in \cite{garholmai07} for the case where
$\xi$ is a symmetric exclusion process in a Bernoulli equilibrium, namely,
between recurrent and transient $p(\cdot,\cdot)$.
\end{enumerate}


\subsection{A scaling conjecture}
\label{S1.9}

Let $p_t(x,y)$ be the probability for the random walk with transition kernel
$p(\cdot,\cdot)$ [satisfying (\ref{pdef})] and step rate 1 to move
from $x$ to
$y$ in time $t$.
The following conjecture is a refinement of Theorem \ref{Lyaprop}(ii)(b).

\begin{conjecture}\label{cj-larg-kappa}
Suppose that $p(\cdot,\cdot)$ is a simple random walk. Then for all
$d\geq5$, $p\in\mathbb{N}$,
$\gamma\in(0,\infty)$ and $\rho\in(0,1)$,
\begin{eqnarray}\label{limlamb*}
&&\lim_{\kappa\to\infty}
2d\kappa[\lambda_p(\kappa)-\rho\gamma]\nonumber\\[-8pt]\\[-8pt]
&&\qquad = \frac{\rho(1-\rho)\gamma^2}{G_d}G_d^\ast
+ 1_{\{d=5\}} (2d)^5 \biggl[\frac{\rho(1-\rho)\gamma^2}{G_d} p \biggr]^2
\mathcal{P}_5\nonumber
\end{eqnarray}
with
%
\begin{equation}\label{GdGdastdefs}
G_d=\int_0^\infty p_t(0,0) \,dt,\qquad
G_d^\ast=\int_0^\infty t p_t(0,0) \,dt
\end{equation}
and
%
\begin{equation}\label{P5def}
\qquad \mathcal{P}_5 = \mathop{\sup_{ f \in H^1(\mathbb{R}^5)}}_{\|f\|_2=1}
\biggl[\int_{\mathbb{R}^5} \int_{\mathbb{R}^5} dx \,dy
\frac{f^2(x) f^2(y)}{16\pi^2\|x-y\|} - \|\nabla f\|_2^2 \biggr] \in
(0,\infty),
\end{equation}
where $\|\cdot\|_2$ is the $L^2$-norm on $\mathbb{R}^5$, $\nabla$ is the
gradient operator,
and $H^1(\mathbb{R}^5)=\{f\dvtx\mathbb{R}^5\to\mathbb{R}\dvtx
f,\nabla f \in L^2(\mathbb{R}^5)\}$.
\end{conjecture}

A remarkable feature of (\ref{limlamb*}) is the occurrence of a
``polaron-type''
term in $d=5$. An important consequence of (\ref{limlamb*}) is that in $d=5$
there exists a $\kappa_1<\infty$ such that
$\lambda_p(\kappa)>\lambda_{p-1}(\kappa)$ for all $\kappa\in
(\kappa_1,\infty)$
when $p=2$ and, by the remark made after formula (\ref{lyasand}), also when
$p\in\mathbb{N}\setminus\{1\}$, that is, the solution of the parabolic
Anderson model is
intermittent for all $\kappa$ sufficiently large. For $d\ge6$,
Conjecture~\ref{cj-larg-kappa} does not allow to decide about intermittency
for large $\kappa$.

The analogue of (\ref{limlamb*}) for independent simple random walks
and simple
symmetric exclusion was proved in \cite{garhol06,garholmai07} and
\cite{garholmai08pr} with quite a bit of effort (with $d=3$ rather
than $d=5$
appearing as the critical dimension). We provide a heuristic
explanation of
(\ref{limlamb*}) in the \hyperref[appA]{Appendix}.


\section[Proof of Theorem 1.3]{Proof of Theorem \protect\ref{Lyaexist}}
\label{S2}

Throughout this section we assume that $p(\cdot,\cdot)$ satisfies
(\ref{pdef}). The
existence of the Lyapunov exponents for $\mu=\nu_\rho S_T$, $T \in
[0,\infty]$,
is proved in Section \ref{S2.1}, the fact that they are equal is
proved in Section
\ref{S2.2}. In what follows, $d\geq1$, $p\in\mathbb{N}$, $\kappa
\in
[0,\infty)$, $\gamma\in
(0,\infty)$ and $\rho\in(0,1)$ are kept fixed. Recall (\ref{lyap2}).

\subsection{Existence of Lyapunov exponents}
\label{S2.1}

\begin{proposition}\label{Lyamunuexist}
For all $T\in[0,\infty]$, the Lyapunov exponent $\lambda^{\nu_\rho S_T}_p$
exists.
\end{proposition}

\begin{pf}
The proof proceeds in 2 steps:

\textit{Step} 1 (Bridge approximation argument).
Let $Q_{t\log t}=\mathbb{Z}^d\cap[-t\log t,t\log t]^d$. As noted in
G\"artner
and den
Hollander \cite{garhol06}, Section 4.1, we have, for $\mu=\nu_\rho S_T$,
%
\begin{equation}\label{supad19}
\qquad \underline{\Lambda}^\mu_p(t)
\leq\Lambda^\mu_p(t)
\leq\frac{1}{pt} \log\bigl( |Q_{t\log t} |^p
e^{pt\underline{\Lambda}^\mu_p(t)}+p e^{\gamma pt}
\mathrm{P}_{0} \bigl(X_1^\kappa(t)\notin Q_{t\log t} \bigr) \bigr)
\end{equation}
with
\begin{eqnarray}\label{supad12}
\underline{\Lambda}^\mu_p(t)
&=&\frac{1}{pt}\log\max_{x\in\mathbb{Z}^d} (\mathbb{E}_{ \mu
}\otimes\mathrm{E}_{
0}^{\otimes p} )\nonumber\\[-8pt]\\[-8pt]
&&{}\times\Biggl(\exp\Biggl[\gamma\int_0^t \sum_{q=1}^{p}\xi
 \bigl(X_q^\kappa(s),t-s \bigr) \,ds \Biggr]
\prod_{q=1}^{p}\delta_x (X_q^\kappa(t) ) \Biggr).\nonumber
\end{eqnarray}
Since $\lim_{t\to\infty}(1/t)\log\mathrm{P}_{0} (X_1^\kappa
(t)\notin
Q_{t\log t} )
=-\infty$, it follows that
%
\begin{equation}\label{supad21}
\lim_{t\rightarrow \infty} [\Lambda^\mu_p(t)-\underline{\Lambda
}^\mu
_p(t) ]=0.
\end{equation}
Hence, to prove the existence of $\lambda^\mu_p$, it suffices to
prove the
existence of
%
\begin{equation}\label{supad23}
\underline{\lambda}^\mu_p=\lim_{t\to\infty}\underline{\Lambda
}^\mu_p(t),
\end{equation}
after which we can conclude from (\ref{supad21})
that $\underline{\lambda}^\mu_p=\lambda^\mu_p$.
We will prove (\ref{supad23})
by showing that $t\mapsto t\underline{\Lambda}^\mu_p(t)$ is superadditive,
which will imply that
%
\begin{equation}\label{supadprop}
\underline{\lambda}^\mu_p = \sup_{t>0} \underline{\Lambda}^\mu_p(t).
\end{equation}

\textit{Step} 2 (Superadditivity).
We first give the proof for $p=1$. To that end, abbreviate
%
\begin{equation}\label{calEdef}
\qquad \mathcal{E}(t,y) = \exp\biggl[\gamma\int_0^t \xi\bigl(X^\kappa(s),t-s \bigr)\, ds \biggr]
\delta_y (X^\kappa(t) ),\qquad  t>0, y\in\mathbb{Z}^d.
\end{equation}
Using formula (\ref{momCRW}) in Proposition \ref{MomCRWprop}, we have,
for all $t_1,t_2>0$ and $x,y\in\mathbb{Z}^d$,
\begin{eqnarray}\label{app1}
&&e^{-\rho\gamma(t_1+t_2)} (\mathbb{E}_{ \nu_\rho S_T} \otimes
\mathrm{E}_{ 0})
\bigl(\mathcal{E}(t_1+t_2,x) \bigr)\nonumber\\
&&\qquad = (\mathrm{E}_{ 0}\otimes\mathbb{E}_{\mathrm{Poiss}}) \bigl(\delta_x
\bigl(X^\kappa(t_1+t_2) \bigr)
\mathbb{E}^\ast\bigl(\rho^{-\mathcal{N}^{\mathrm{coal}}_{T+t_1+t_2} \{
(X^\kappa(s),s)\dvtx s\in\Pi_{\rho\gamma}([0,t_1+t_2]) \}} \bigr) \bigr)\nonumber\\
&&\qquad \geq(\mathrm{E}_{ 0}\otimes\mathbb{E}_{\mathrm{Poiss}}) \bigl(\delta_y
(X^\kappa(t_1) )
\delta_x \bigl(X^\kappa(t_1+t_2) \bigr)\nonumber\\
&&\hphantom{(\mathrm{E}_{ 0}\otimes\mathbb{E}_{\mathrm{Poiss}}) \bigl(}\qquad \quad {} \times
\mathbb{E}^\ast\bigl(
\rho^{-\mathcal{N}^{\mathrm{coal}}_{T+t_1} \{(X^\kappa(s),s)\dvtx
s\in\Pi_{\rho
\gamma}([0,t_1]) \}}\nonumber\\[-8pt]\\[-8pt]
&&\hphantom{(\mathrm{E}_{ 0}\otimes\mathbb{E}_{\mathrm{Poiss}}) \bigl({} \times\mathbb{E}^\ast\bigl(}
\qquad \quad {} \times
\rho^{-\mathcal{N}^{\mathrm{coal}}_{T+t_1+t_2} \{(X^\kappa
(s),s)\dvtx
s\in\Pi_{\rho\gamma}([t_1,t_1+t_2]) \}} \bigr) \bigr)\nonumber\\
&&\qquad = (\mathrm{E}_{ 0}\otimes\mathbb{E}_{\mathrm{Poiss}}) \bigl(\delta_y
(X^\kappa(t_1) )
\delta_{x-y} \bigl(X^\kappa(t_1+t_2)-X^\kappa(t_1) \bigr)\nonumber\\
&&\hphantom{(\mathrm{E}_{ 0}\otimes\mathbb{E}_{\mathrm{Poiss}}) \bigl(}\qquad \quad {} \times
\mathbb{E}^\ast\bigl(\rho^{-\mathcal{N}^{\mathrm{coal}}_{T+t_1} \{
(X^\kappa(s),s)\dvtx s\in\Pi_{\rho\gamma}([0,t_1]) \}}\nonumber\\
&&\hphantom{(\mathrm{E}_{ 0}\otimes\mathbb{E}_{\mathrm{Poiss}}) \bigl({} \times\mathbb{E}^\ast\bigl(}
\qquad \quad {} \times
\rho^{-\mathcal{N}^{\mathrm{coal}}_{T+t_1+t_2} \{(X^\kappa
(s)-X^\kappa(t_1),s)\dvtx
s\in\Pi_{\rho\gamma}([t_1,t_1+t_2]) \}} \bigr) \bigr),\nonumber
\end{eqnarray}
where the inequality comes from inserting the extra factor $\delta
_y(X^\kappa(t_1))$
under the expectation and ignoring coalescence between random walks
that start before,
respectively, after time $t_1$, and the last line uses the
shift-invariance of
$\mathcal{N}^{\mathrm{coal}}_{T+t_1+t_2}$. Because $X^\kappa$ and $\Pi
_{\rho\gamma
}$ have
independent stationary increments, we have
\begin{eqnarray}\label{app3}
&&\mbox{r.h.s. (\ref{app1})}\nonumber\\
&&\qquad = (\mathrm{E}_{ 0}\otimes\mathbb{E}_{\mathrm{Poiss}}) \bigl(\delta_y
(X^\kappa(t_1) )
\mathbb{E}^\ast\bigl(\rho^{-\mathcal{N}^{\mathrm{coal}}_{T+t_1} \{
(X^\kappa(s),s)\dvtx
s\in\Pi_{\rho\gamma}([0,t_1]) \}} \bigr)
\bigr)\nonumber\\[-8pt]\\[-8pt]
&&\qquad\quad {}  \times(\mathrm{E}_{ 0}\otimes\mathbb{E}_{\mathrm{Poiss}}) \bigl(\delta
_{x-y} (X^\kappa(t_2) )
\mathbb{E}^\ast\bigl(\rho^{-\mathcal{N}^{\mathrm{coal}}_{T+t_2} \{
(X^\kappa(s),s)\dvtx
s\in\Pi_{\rho\gamma}([0,t_2]) \}} \bigr) \bigr)\nonumber\\
&&\qquad = e^{-\rho\gamma t_1} (\mathbb{E}_{ \nu_\rho S_T} \otimes
\mathrm{E}_{ 0}) (\mathcal{E}(y,t_1) ) \times e^{-\rho\gamma t_2}
(\mathbb{E}_{ \nu_\rho S_T}
\otimes\mathrm{E}_{ 0}) \bigl(\mathcal{E}(x-y,t_2) \bigr),\nonumber
\end{eqnarray}
where in the last line we again use formula (\ref{momCRW}).
Taking the maximum over $x,y\in\mathbb{Z}^d$ in (\ref{app1})--(\ref{app3}),
we conclude that
%
\begin{equation}\label{app12}
\qquad \exp[(t_1+t_2)\underline{\Lambda}^{\nu_\rho S_T}_1(t_1+t_2)]
\ge\exp[t_1\underline{\Lambda}^{\nu_\rho S_T }_1(t_1)] \times
\exp[t_2\underline{\Lambda}^{\nu_\rho S_T}_1(t_2)],
\end{equation}
which proves the superadditivity of $t\mapsto
t\underline{\Lambda}^{\nu_\rho S_T}_1(t)$.

The same proof works for $p\in\mathbb{N}\setminus\{1\}$. Simply replace
(\ref{calEdef}) by
%
\begin{eqnarray}\label{calEdefp}
\mathcal{E}_p(t,y) = \exp\Biggl[\gamma\int_0^t \sum_{q=1}^p
\xi\bigl(X^\kappa_q(s),t-s \bigr) \,ds \Biggr]
\prod_{q=1}^p \delta_y (X^\kappa_q(t) ),\nonumber\\[-8pt]\\[-8pt]
\eqntext{t \geq0, y\in\mathbb{Z}^d,}
\end{eqnarray}
and proceed in a similar manner.
\end{pf}


\subsection{Equality of Lyapunov exponents}
\label{S2.2}

\begin{proposition}\label{lpmunuub}
$\lambda_p^{\nu_\rho} = \lambda_p^{\nu_\rho S_T}$ for all $T\in
[0,\infty]$.
In particular, $\lambda^{\nu_\rho}=\lambda^{\mu_\rho}$.
\end{proposition}

\begin{pf}
We first give the proof for $p=1$.

$\lambda_1^{\nu_\rho} \leq\lambda_1^{\nu_\rho S_T}$:
Since $t \mapsto\mathcal{N}^{\mathrm{coal}}_t$ is nondecreasing, it is
immediate from the
representation formula (\ref{lyaCRW}) in Proposition \ref{MomCRWprop} that
%
\begin{equation}\label{Lambdamunuub}
\Lambda_1^{\nu_\rho}(t) \leq\Lambda_1^{\nu_\rho S_T}(t)
\qquad \forall t>0, T\in[0,\infty].
\end{equation}
Since $\lambda_1^{\nu_\rho S_T}=\lim_{t\to\infty}\Lambda_1^{\nu
_\rho S_T}(t)$,
this implies the claim.

$\lambda_1^{\nu_\rho} \geq\lambda_1^{\nu_\rho S_T}$:
We first assume that $T<\infty$.
Recall (\ref{supad21}) and (\ref{supad23})--(\ref{calEdef}), and estimate,
for $T,t>0$,
\begin{eqnarray}\label{lbmon1}
\lambda_1^{\nu_\rho} &=& \Lambda_1^{\nu_\rho}(\infty)
= \underline{\Lambda}_1^{ \nu_\rho}(\infty)\nonumber\\[-8pt]\\[-8pt]
&\geq&\underline{\Lambda}_1^{\nu_\rho}(T+t)
= \frac{1}{T+t} \log\max_{x\in\mathbb{Z}^d}
(\mathbb{E}_{ \nu_\rho} \otimes\mathrm{E}_{ 0}) \bigl(\mathcal
{E}(T+t,x) \bigr).\nonumber
\end{eqnarray}
In the right-hand side of (\ref{lbmon1}), drop the part $s \in[t,T+t]$
from the integral over $s \in[0,T+t]$ in definition (\ref{calEdef}) of
$\mathcal{E}(T+t,x)$, insert an extra factor $\delta_x(X^\kappa(t))$ under
the expectation,
and use the Markov property of $\xi$ and $X^\kappa$ at time $t$. This gives
%
\begin{equation}\label{lammon2}
\quad\qquad \mbox{r.h.s. (\ref{lbmon1})} \geq\frac{1}{T+t} \log\max_{x\in
\mathbb{Z}^d}
\bigl\{
(\mathbb{E}_{ \nu_\rho S_T} \otimes\mathrm{E}_{ 0} ) (\mathcal
{E}(t,x) )
\mathrm{P}_0 \bigl(X^\kappa(T)=0 \bigr) \bigr\}.
\end{equation}
Combine (\ref{lbmon1}) with (\ref{lammon2}) to get
%
\begin{equation}\label{lammon3}
\lambda_1^{\nu_\rho} \geq
\frac{t}{T+t} \underline{\Lambda}_1^{\nu_\rho S_T}(t)
+ \frac{1}{T+t} \log\mathrm{P}_0 \bigl(X^\kappa(T)=0 \bigr).
\end{equation}
Let $t\to\infty$ to get $\lambda_1^{\nu_\rho} \geq\underline
{\Lambda}_1^{\nu_\rho S_T}
(\infty)= \lambda_1^{\nu_\rho S_T}$, which proves the claim.

Next, for $T,t>0$ and $x\in\mathbb{Z}^d$,
%
\begin{equation}\label{lammon4}
\qquad \lambda_1^{\nu_\rho}\ge\lambda_1^{\nu_\rho S_T}
= \underline{\Lambda}_1^{\nu_\rho S_T}(\infty)
\geq\underline{\Lambda}_1^{\nu_\rho S_T}(t)
\geq\frac{1}{t} \log(\mathbb{E}_{ \nu_\rho S_T} \otimes\mathrm
{E}_{ 0} ) (\mathcal{E}
(t,x) ),
\end{equation}
where we have used (\ref{supadprop}). The weak convergence of $\nu
_\rho S_T$ to $\mu_\rho$ implies that we can
take the limit as $T\to\infty$ to obtain
%
\begin{equation}\label{lammon6}
\lambda_1^{\nu_\rho} \geq\frac{1}{t}
\log(\mathbb{E}_{ \mu_\rho}\otimes\mathrm{E}_{ 0} ) (\mathcal
{E}(t,x) ).
\end{equation}
Finally, taking the maximum over $x$ and letting $t\to\infty$, we arrive
at $\lambda_1^{\nu_\rho}\ge\lambda_1^{\mu_\rho}$, which is the
claim for
$T=\infty$.

The same proof works for $p\in\mathbb{N}\setminus\{1\}$ by using
(\ref{calEdefp})
instead of (\ref{calEdef}).
\end{pf}


\section[Proof of Theorem 1.4]{Proof of Theorem \protect\ref{Lyalip}}
\label{S3}

Throughout this section we assume that $p(\cdot,\cdot)$ satisfies
(\ref{pdef}).
In Section \ref{S3.1} we show that $\kappa\mapsto\lambda_p(\kappa
)$ is globally
Lipschitz outside any neighborhood of $0$. In Section \ref{S3.2} we
show that
$\lambda_p(\kappa)>\rho\gamma$ for all $\kappa\in[0,\infty)$. In
what follows,
$d \geq1$, $p\in\mathbb{N}$, $\gamma\in(0,\infty)$ and $\rho\in(0,1)$
are kept
fixed.

\subsection{Lipschitz continuity}
\label{S3.1}

In this section we prove Theorem \ref{Lyalip}(i).

\begin{pf*}{Proof of Theorem \ref{Lyalip}\textup{(i)}}
In what follows, $\mu$ can be any of the initial distributions $\nu
_\rho S_T$,
$T\in[0,\infty]$ (recall Proposition \ref{lpmunuub}).
We write $\Lambda_p^\mu(\kappa;t)$ to indicate the
$\kappa$-dependence
of $\Lambda_p^\mu(t)$ given by (\ref{fey-kac2}). We give the proof
for $p=1$.

Pick $\kappa_1,\kappa_2\in(0,\infty)$ with $\kappa_1<\kappa_2$
arbitrarily.
By a standard application of Girsanov's formula,
%
\begin{eqnarray}\label{Gir}
&&\exp[t\Lambda^\mu_1(\kappa_2;t)]\nonumber\\
&&\qquad  = (\mathbb{E}_{ \mu} \otimes\mathrm{E}_{ 0})
\biggl(\exp\biggl[\gamma\int_0^t \xi\bigl(X^{\kappa_2}(s),t-s\bigr) \,ds \biggr] \biggr)\nonumber\\
&&\qquad  = (\mathbb{E}_{ \mu} \otimes\mathrm{E}_{ 0})
\biggl(\exp\biggl[\gamma\int_0^t \xi\bigl(X^{\kappa_1}(s),t-s\bigr) \,ds \biggr] \\
&&\hphantom{(\mathbb{E}_{ \mu} \otimes\mathrm{E}_{ 0})\biggl(}\qquad \quad {} \times
\exp[J(X^{\kappa_1};t)\log(\kappa_2/\kappa_1)-2d(\kappa_2-\kappa
_1)t ] \biggr)\nonumber\\
&&\qquad  = I+\mathit{II},\nonumber
\end{eqnarray}
where $J(X^{\kappa_1};t)$ is the number of jumps of $X^{\kappa_1}$ up
to time $t$,
$I$ and $\mathit{II}$ are the contributions coming from the events $\{
J(X^{\kappa_1};t)\leq
M2d\kappa_2t\}$, respectively, $\{J(X^{\kappa_1};t)> M2d\kappa_2t\}
$, and $M>1$ is
to be chosen. Clearly,
%
\begin{equation}\label{Iest}
I \leq\exp\bigl[ \bigl(M2d\kappa_2\log(\kappa_2/\kappa_1)
-2d(\kappa_2-\kappa_1) \bigr)t \bigr]
\exp[t\Lambda^\mu_1(\kappa_1;t)],
\end{equation}
while
%
\begin{equation}\label{IIest}
\mathit{II} \leq e^{\gamma t} \mathrm{P}_{0} \bigl(J(X^{\kappa_2};t)>M2d\kappa_2t \bigr)
\end{equation}
because we may estimate $\int_0^t \xi(X^{\kappa_1}(s),t-s) \,ds \leq t$
and afterward use Girsanov's formula in the reverse direction.
Since $J(X^{\kappa_2};t)=J^*(2d\kappa_2t)$
with $(J^*(t))_{t\geq0}$ a rate-$1$ Poisson process, we have
%
\begin{equation}\label{LDPPoi}
\lim_{t\to\infty} \frac{1}{t} \log\mathrm{P}_{0} \bigl(J(X^{\kappa
_2};t)>M2d\kappa_2t \bigr)
= -2d\kappa_2 \mathcal{I}(M)
\end{equation}
with
%
\begin{equation}\label{IMid}
\mathcal{I}(M) = \sup_{u\in\mathbb{R}} [Mu- (e^u-1 ) ] = M\log M-M+1.
\end{equation}
Since $\lambda_1(\kappa)=\lim_{t\to\infty} \Lambda_1^\mu(\kappa
;t)$, it follows
from (\ref{Gir})--(\ref{LDPPoi}) that
\begin{eqnarray}\label{lk2lk1bd}
\lambda_1(\kappa_2) &\leq&
[M2d\kappa_2\log(\kappa_2/\kappa_1)
-2d(\kappa_2-\kappa_1)+\lambda_1(\kappa_1) ]\nonumber\\[-8pt]\\[-8pt]
&&{}\vee[\gamma-2d\kappa_2\mathcal{I}(M) ].\nonumber
\end{eqnarray}
On the other hand, estimating $J(X^{\kappa_1};t)\geq0$ in (\ref{Gir}),
we have
%
\begin{equation}\label{Jt0}
\exp[t\Lambda^\mu_1(\kappa_2;t)]
\geq\exp[-2d(\kappa_2-\kappa_1)t] \exp[t\Lambda^\mu_1(\kappa_1;t)],
\end{equation}
which gives the lower bound
%
\begin{equation}\label{lk2lk1lb}
\lambda_1(\kappa_2) - \lambda_1(\kappa_1) \geq-2d(\kappa_2-\kappa_1).
\end{equation}

Next, for $\kappa\in(0,\infty)$, define
\begin{eqnarray}\label{derdefs}
D^+\lambda_1(\kappa) &=& \limsup_{\delta\to0}
\delta^{-1}
[\lambda_1(\kappa+\delta)-\lambda_1(\kappa)],\nonumber\\[-8pt]\\[-8pt]
D^-\lambda_1(\kappa) &=& \liminf_{\delta\to0}
\delta^{-1} [\lambda_1(\kappa+\delta)-\lambda_1(\kappa)].\nonumber
\end{eqnarray}
Then, picking $\kappa_1=\kappa$ and $\kappa_2=\kappa+\delta$
(resp.,
$\kappa_1=\kappa-\delta$ and $\kappa_2=\kappa$) in (\ref{lk2lk1bd})
and letting $\delta\downarrow0$, we get
%
\begin{equation}\label{D1+ub}
D^+\lambda_1(\kappa) \leq(M-1)2d \qquad \forall M>1\dvtx
2d\kappa\mathcal{I}(M)-(1-\rho)\gamma\geq0
\end{equation}
[with the latter together with $\lambda_1(\kappa)\ge\rho\gamma$
guaranteeing
that the first term in the right-hand side of (\ref{lk2lk1bd}) is
the maximum], while (\ref{lk2lk1lb}) gives
%
\begin{equation}\label{D1-lb}
D^-\lambda_1(\kappa) \geq-2d.
\end{equation}
We may pick
%
\begin{equation}\label{Mcond}
M = M(\kappa) = \mathcal{I}^{-1} \biggl(\frac{(1-\rho)\gamma}{2d\kappa} \biggr)
\end{equation}
with $\mathcal{I}^{-1}$ the inverse of $\mathcal{I}\dvtx[1,\infty
)\to\mathbb{R}$.
Since $\mathcal{I}(M)=\frac12(M-1)^2[1+o(1)]$ as $M \downarrow1$, it
follows that
%
\begin{equation}\label{D+asymp}
[M(\kappa)-1]2d = 2d \sqrt{\gamma\frac{1-\rho}{d\kappa}} [1+o(1)]
\qquad \mbox{as } \kappa\to\infty.
\end{equation}
By (\ref{D1+ub}), the latter implies that $\kappa\mapsto D^+\lambda
_1(\kappa)$ is
bounded from above outside any neighborhood of $0$. Since, by (\ref
{D1-lb}), $\kappa\mapsto
D^-\lambda_1(\kappa)$ is bounded from below, the claim follows.

The extension to $p\in\mathbb{N}\setminus\{1\}$ is straightforward
and is
left to the reader.
\end{pf*}


\subsection{Clumping}
\label{S3.2}

In this section we prove Theorem \ref{Lyalip}(ii).

\begin{pf*}{Proof of Theorem \ref{Lyalip}\textup{(ii)}}
Fix $d \geq1$, $\kappa\in[0,\infty)$, $\gamma\in(0,\infty)$ and
$\rho\in(0,1)$. Since $p\mapsto\lambda_p(\kappa)$ is
nondecreasing, it
suffices to give the proof for $p=1$. In what follows, $\mu$ can be
any of
the measures $\nu_\rho S_T$, $T\in[0,\infty]$
(recall Proposition~\ref{lpmunuub}).

Abbreviate
%
\begin{equation}\label{Iabbrv}
I(X^\kappa;T) = \gamma\int_0^T ds\, \bigl[\xi\bigl(X^\kappa(s),T-s \bigr)-\rho\bigr],
\qquad T>0.
\end{equation}
For any $T>0$ we have, recalling (\ref{supad12})--(\ref{supadprop}),
\begin{eqnarray}\label{lambda-lb9}
\lambda_1(\kappa) &=& \Lambda^\mu_1(\infty)
= \underline{\Lambda}^\mu_1(\infty) \geq\underline{\Lambda}^\mu
_1(T)\nonumber\\
&\geq&\rho\gamma+ \frac1T \log(\mathbb{E}_{ \mu}\otimes\mathrm
{E}_{ 0}) (\exp
[I(X^\kappa;T)]
\delta_0 (X^\kappa(T) ) )\nonumber\\[-8pt]\\[-8pt]
&\geq&\rho\gamma+ \frac1T \log(\mathbb{E}_{ \mu}\otimes\mathrm
{E}_{ 0}) \biggl( \biggl[1
+ I(X^\kappa;T)\nonumber\\
&&\hphantom{\rho\gamma+ \frac1T \log(\mathbb{E}_{ \mu}\otimes\mathrm
{E}_{ 0}) \biggl( \biggl[}
{}+ \frac12 I(X^\kappa;T)^2e^{-\gamma T} \biggr] \delta_0 (X^\kappa(T) ) \biggr),\nonumber
\end{eqnarray}
where in the third line we use that $e^x \geq1+x+\frac12
x^2e^{-|x|}$, $x\in\mathbb{R}$.

As $T\downarrow0$, we have
\begin{eqnarray}\label{lambda-lb11}
\qquad (\mathbb{E}_{ \mu}\otimes\mathrm{E}_{ 0}) \biggl( \biggl[\frac1T I(X^\kappa
;T) \biggr]^2
\delta_0 (X^\kappa(T) ) \biggr) &\to&
\gamma^2\int_{\Omega}\mu(d\eta)
[\eta(0)-\rho]^2\nonumber\\[-8pt]\\[-8pt]
&=& \rho(1-\rho)\gamma^2\nonumber
\end{eqnarray}
and
%
\begin{equation}\label{lambda-lb11*}
(\mathbb{E}_{ \mu}\otimes\mathrm{E}_{ 0}) \biggl( \biggl[\frac1T I(X^\kappa
;T) \biggr]
\delta_0 (X^\kappa(T) ) \biggr)
\ge-O(T^2).
\end{equation}
The claim in (\ref{lambda-lb11}) is obvious, the claim in (\ref{lambda-lb11*}) will be proven
below. Combining (\ref{lambda-lb9})--(\ref{lambda-lb11*}), we have
%
\begin{equation}\label{lambda-lb12}
\lambda_1(\kappa)-\rho\gamma\geq\tfrac14T \rho(1-\rho)\gamma^2,
\qquad 0 < T \leq T_0(\kappa),
\end{equation}
for some $T_{0}(\kappa)<\infty$, showing that $\lambda_1(\kappa
)>\rho\gamma$.

To prove (\ref{lambda-lb11*}), let $J(X^\kappa;T)$ denote the number
of jumps
by $X^\kappa$ up to time $T$. Then
%
\begin{eqnarray}\label{lambda-lb13}
&&(\mathbb{E}_{ \mu}\otimes\mathrm{E}_{ 0}) \biggl( \biggl[\frac1T I(X^\kappa
;T) \biggr]
\delta_0 (X^\kappa(T) ) \biggr)\nonumber\\
&&\qquad  = (\mathbb{E}_{ \mu}\otimes\mathrm{E}_{ 0}) \biggl( \biggl[\frac1T
I(X^\kappa;T) \biggr]
\delta_0 (X^\kappa(T) ) \\
&&\hphantom{(\mathbb{E}_{ \mu}\otimes\mathrm{E}_{ 0}) \biggl(}\qquad \quad {} \times
\bigl(\mathbh{1}\{J(X^\kappa;T)=0\} + \mathbh{1}\{J(X^\kappa;T) \geq1\}
\bigr) \biggr).\nonumber
\end{eqnarray}
The first term in the right-hand side of (\ref{lambda-lb13}) equals
%
\begin{equation}\label{lambda-lb14}
\mathrm{P}_0 \bigl(J(X^\kappa;T)=0 \bigr) \frac{\gamma}{T} \int_0^T
ds\,
\mathbb{E}_{ \mu} \bigl(\xi(0,s)-\rho\bigr) = 0,
\end{equation}
while the second term is bounded below by
\begin{eqnarray}\label{lambda-lb15}
-\rho\gamma\mathrm{P}_0 \bigl(J(X^\kappa;T) \geq1, X^\kappa(T)=0 \bigr)
&\geq&-\rho\gamma\mathrm{P}_0 \bigl(J(X^\kappa;T) \geq2 \bigr)\nonumber\\[-8pt]\\[-8pt]
&=& -O(T^2),\nonumber
\end{eqnarray}
as $T \downarrow0$. Combine (\ref{lambda-lb13})--(\ref{lambda-lb15})
to get the
claim in (\ref{lambda-lb11*}).
\end{pf*}


\section[Proof of Theorem 1.5(ii)(a) and (b)]{Proof of Theorem \protect\ref{Lyaprop}(ii)(a) and (b)}
\label{S4}

Throughout this section we assume that $p(\cdot,\cdot)$ satisfies
(\ref{pdef})
and that $d \geq5$. In Section \ref{S4.1} we state an estimate for
blocks of
coalescing random walks. In Section \ref{S4.2} we formulate two
lemmas, and in
Section \ref{S4.3} we use these lemmas to prove the block estimate.
The block
estimate is used in Sections \ref{S4.4} and \ref{S4.5} to prove Theorem
\ref{Lyaprop}(ii)(a) and (b), respectively.


\subsection{Block estimate}
\label{S4.1}

We call a collection of subsets $S_1,\ldots,S_N$ of $\mathbb{R}$
\textit{ordered},
if $s<t$ for all $s\in S_i$, $t\in S_j$ and $i<j$.
Given a path $\psi\dvtx\mathbb{R}\to\mathbb{Z}^d$ and a
collection of disjoint
finite subsets $S_1,\ldots,S_N$ of $\mathbb{R}$,
we are going to estimate the moment generating function of
$\mathcal{N}^{\mathrm{coal}}_\infty\{(\psi(s),s)\dvtx s\in
\bigcup_{j=1}^N S_j \}$, the number of random walks starting from sites
$\psi(s)$ at times $s\in\bigcup_{j=1}^N S_j$ that coalesce eventually\vspace*{1pt}
[recall (\ref{Ntcoaldef})].
Let $d(S_i,S_j)$ denote the Euclidean distance between $S_i$ and $S_j$.

Our key estimate, which will be proved in Section \ref{S4.3}, is the following proposition.

\begin{proposition}\label{blockprop}
Let $d \geq5$. Then there exist $\delta\dvtx(0,\infty)\to
(0,\infty)$ with
$\lim_{K\to\infty}\delta(K)=0$ and, for each $\epsilon\in(0,(d-4)/2)$,
$C_\epsilon>0$ such that the following holds.
For all $\rho\in(0,1)$, $\psi\dvtx\mathbb{R}\to\mathbb{Z}^d$,
all ordered collections
of disjoint finite subsets $S_1,\ldots,S_N$ of $\mathbb{R}$, all
$\epsilon\in(0,(d-4)/2)$, $K>0$ and $r,r'>1$ with $1/r+1/r'=1$,
%
\begin{eqnarray}\label{block}
&&\mathbb{E}^\ast\bigl(
\rho^{-\mathcal{N}^{\mathrm{coal}}_\infty\{(\psi(s),s)\dvtx s\in
\bigcup_{j=1}^N S_j \}} \bigr)\nonumber\\
&&\qquad  \leq\exp\Biggl[\frac{\delta(K)}{\rho}\sum_{j=1}^N |S_j|
+ C_\epsilon K \frac{\rho^{-r'}-1}{r'}
\sum_{1\le j<k \le N} \frac{|S_j| |S_k|}{d(S_j,S_k)^{1+\epsilon}} \Biggr]\\
&&\qquad \quad {} \times
\Biggl[ \prod_{j=1}^N \mathbb{E}^\ast\bigl(\rho^{-r \mathcal{N}^\mathrm
{coal}_\infty\{(\psi(s),s)\dvtx s\in S_j \}} \bigr)\Biggr]^{1/r}.\nonumber
\end{eqnarray}
\end{proposition}

Let $I_1^\prime,I_1^{\prime\prime},\ldots,I_N^\prime,I_N^{\prime
\prime}$ be
a finite collection of adjacent time intervals and assume that
$S_j\subset I_j^\prime$ for $j=1,\ldots,N$. What the above
proposition does
is \textit{decouple} the coalescing random walks that start in disjoint
\textit{time-blocks} $I_j^\prime$ separated by \textit{time-gaps}
$I_j^{\prime\prime}$.


\subsection{Preparatory lemmas}
\label{S4.2}

To prove Proposition \ref{blockprop}, we need Lemmas \ref{badlm}--\ref{Ycondlem}
below. To this end, fix a path $\psi\dvtx\mathbb{R}\to\mathbb
{Z}^d$ arbitrarily.
Let $(Y^u)_{u \in\mathbb{R}}$ be a family of independent random walks
$Y^u$ with
transition kernel $p^\ast(\cdot,\cdot)$ and step rate $1$ starting from
$\psi(u)$ at time $u$. Set $Y^u(s)=\psi(u)$ for $s<u$. We write
$\mathbb{P}
^*$ for the
joint law of these random walks.

Given $u\in\mathbb{R}$ and $j\in\mathbb{Z}$, let
%
\begin{equation}\label{defrange}
R^u_j = \{Y^u(s)\dvtx s\in[j,j+1] \}
\end{equation}
denote the \textit{range} of $Y^u$ in the time interval $[j,j+1]$.
For $u\in\mathbb{R}$ and $K>0$, define the event that $Y^u$ is \emph
{$K$-good}
by
%
\begin{equation}\label{Gudef}
G_K^u = \bigcap_{j=\lfloor u\rfloor}^{\infty}
\{ |R^u_j | \le K\log(j-\lfloor u\rfloor+5) \}.
\end{equation}
For $u,v\in\mathbb{R}$ with $u<v$, define the event that $Y^u$ and $Y^v$
\textit{meet} by
%
\begin{equation}\label{Ymeet}
M^{u,v} = \{\exists s\geq v\dvtx Y^u(s)=Y^v(s) \}.
\end{equation}

Our two lemmas stated below give bounds for the probabilities of random walks
not to be $K$-good, respectively, to meet given that the random walk that
starts later is $K$-good.

\begin{lemma}\label{badlm}
For all $u\in\mathbb{R}$ and $K>0$,
%
\begin{equation}\label{badeq1}
\mathbb{P}^\ast([G_K^u]^c)\leq\delta(K)
\end{equation}
with
%
\begin{equation}\label{badeq2}
\delta(K)=\sum_{j=5}^\infty
\exp\bigl[-\lfloor K\log j\rfloor
\bigl(\log(\lfloor K\log j\rfloor)-1 \bigr)-1 \bigr]<\infty
\end{equation}
satisfying $\lim_{K\to\infty}\delta(K)=0$.
\end{lemma}

\begin{pf}
Recalling (\ref{Gudef}) and taking into account that $Y^u$ has stationary
increments, we have
%
\begin{equation}\label{badeq5}
\quad \mathbb{P}^\ast([G_K^u]^c)
\leq\sum_{j=0}^{\infty}
\mathbb{P}^\ast\bigl( |R^0_j | > K\log(j+5) \bigr)
\leq\sum_{j=5}^{\infty}\mathbb{P}^\ast(N_1\geq\lfloor K\log
j\rfloor),
\end{equation}
where $N_1$ denotes the Poisson number of jumps of $Y^0$ during a time interval
of length $1$. An application of Chebyshev's exponential inequality yields,
for $\beta>0$,
%
\begin{eqnarray}\label{badeq7}
\mathbb{P}^\ast(N_1\geq\lfloor K\log j\rfloor)
&\leq& e^{-\beta\lfloor K\log j\rfloor} \mathbb{E}^\ast(e^{\beta
N_1} )\nonumber\\
&=&\exp[-\beta\lfloor K\log j\rfloor+e^\beta-1 ]\\
&=&\exp\bigl[-\lfloor K\log j\rfloor\bigl(\log(\lfloor K\log j\rfloor)-1 \bigr)-1 \bigr],\nonumber
\end{eqnarray}
where in the last line we optimize over the choice of $\beta$ by taking
$\beta=\log(\lfloor K\times\log j\rfloor)$.
Combining (\ref{badeq5}) and (\ref{badeq7}), we get the claim.
\end{pf}

\begin{lemma}\label{Ycondlem}
Let $d \geq5$. Then for all $\epsilon\in(0,(d-4)/2)$ there exists
$C_\epsilon>0$ such that for all $K>0$ and all $u,v\in\mathbb{R}$
with $u<v$,
%
\begin{equation}\label{P*Cest}
\mathbb{P}^\ast(M^{u,v} \mid Y^v )
\le\frac{C_\epsilon K}{(v-u)^{1+\epsilon}}\qquad
\mbox{on }G_K^v.
\end{equation}
\end{lemma}

\begin{pf}
Fix $u,v\in\mathbb{R}$ with $u<v$. Recall (\ref{defrange})--(\ref{Ymeet})
to see that
%
\begin{equation}\label{Cuvinc}
M^{u,v}\subseteq\bigcup_{j=\lfloor v \rfloor}^\infty
\bigcup_{z\in R^v_j} \{\exists s\in[j,j+1]\dvtx Y^u(s)=z \}.
\end{equation}
Hence,
%
\begin{equation}\label{Cuv1}
\mathbb{P}^\ast(M^{u,v}\mid Y^v )
\le\sum_{j=\lfloor v \rfloor}^\infty\sum_{z\in R^v_j}
\mathbb{P}^\ast\bigl(\exists s\in[j,j+1]\dvtx Y^u(s)=z \bigr).
\end{equation}
Since the transition kernel $p^\ast(\cdot,\cdot)$ generates $\mathbb{Z}^d$
[recall (\ref{pdef})], there exists a constant $C>0$ such that
%
\begin{equation}\label{Spitzest}
p^\ast_t(x,y) \leq\frac{C}{(t+7)^{d/2}}
\qquad \forall t\ge0,\  \forall x,y\in\mathbb{Z}^d
\end{equation}
(see Spitzer \cite{sp76}, Proposition 7.6). Let $Y$ be a random walk
on $\mathbb{Z}^d$ with
transition kernel $p^\ast(\cdot,\cdot)$ and jump rate $1$. Let
$\mathbb{P}
^Y_y$ denote its
law when starting at $y$ and $\tau_z=\inf\{s\ge0\dvtx Y(s)=z\}$ its
first hitting
time of $z$. Then, since $Y^u$ and $Y$ have the same independent and
stationary increments,
we have, for $j\ge\lfloor v\rfloor$,
%
\begin{eqnarray}\label{Cuv2}
\mathbb{P}^\ast\bigl(\exists s\in[j,j+1]\dvtx Y^u(s)=z \bigr)
&\le&\sum_{y\in\mathbb{Z}^d} p^\ast_{(j\vee u)-u}(\psi(u),y)
\mathbb{P}^Y_y (\tau
_z\le1 )\nonumber\\
&\le&\frac{C}{(j-u+6)^{d/2}} \sum_{y\in\mathbb{Z}^d} \mathbb
{P}^Y_0 (\tau_y\le1 )\\
&=& \frac{C}{(j-u+6)^{d/2}} \mathbb{E}^Y_0(|R|),\nonumber
\end{eqnarray}
where $R=\{Y(s)\dvtx s\in[0,1]\}$ is the range of $Y$ in the time
interval $[0,1]$.
Since $|R|\leq1+N_1$ with $N_1$ the Poisson number of jumps of $Y$ in
$[0,1]$, we
have $\mathbb{E}^Y_0(|R|) \le2$. Now assume that $Y^v$ is $K$-good [recall
(\ref{Gudef})].
Then, combining (\ref{Cuv1}) with (\ref{Cuv2}),
we obtain
%
\begin{eqnarray}\label{Cuv3}
\mathbb{P}^\ast(M^{u,v}\mid Y^v ) &\le&2C K \sum_{j=\lfloor v
\rfloor
}^\infty
\frac{\log(j-\lfloor v \rfloor+5)}{(j-u+6)^{d/2}}\nonumber\\
&\le&2C K \sum_{j=\lfloor v \rfloor}^\infty
\frac{\log(j-\lfloor u \rfloor+5)}{(j-\lfloor u \rfloor+5)^{d/2}}\\
&\le&2C K \frac{\log(\lfloor v\rfloor-\lfloor u\rfloor+4)}
{(\lfloor v\rfloor-\lfloor u\rfloor+4)^{(d-2)/2}}.\nonumber
\end{eqnarray}
Since $d\ge5$, this clearly implies (\ref{P*Cest}).
\end{pf}


\subsection{Proof of block estimate}
\label{S4.3}

In this section we use Lemmas \ref{badlm} and \ref{Ycondlem} to prove Proposition \ref{blockprop}.

\begin{pf*}{Proof of Proposition \ref{blockprop}}
Fix a path $\psi\dvtx\mathbb{R}\to\mathbb{Z}^d$ and an ordered
collection of
disjoint finite
subsets $S_1,\ldots,S_N$ of $\mathbb{R}$ arbitrarily. Assume that the
coalescing
random walks starting from sites $\psi(s)$ at times $s\in\bigcup
_{j=1}^N S_j$
are constructed from the independent random walks
$Y^u$, $u\in\bigcup_{j=1}^N S_j$, introduced in Section \ref{S4.2},
in the obvious recursive manner: if two walks meet for the first time,
then the random walk that started earlier is killed and the random walk that
started later survives.

Now recall (\ref{Gudef}). Distinguishing between all possible ways to
distribute the good and the bad events and using the independence of the
random walks $Y^u$, we estimate
%
\begin{eqnarray}\label{block1}
&&\mathbb{E}^\ast\bigl( \rho^{-\mathcal{N}^{\mathrm{coal}}_\infty\{(\psi
(s),s)\dvtx
s\in\bigcup_{j=1}^N S_j \}} \bigr)\nonumber\\
&&\qquad  = \mathop{\sum_{A_i\subseteq S_i} }_{1\leq i\leq N}
\mathbb{E}^\ast\Biggl(\rho^{-\mathcal{N}^{\mathrm{coal}}_\infty\{(\psi
(s),s)\dvtx
s\in\bigcup_{j=1}^N S_j \}}\nonumber \\
&&\qquad\quad\hphantom{\mathop{\sum_{A_i\subseteq S_i} }_{1\leq i\leq N}
\mathbb{E}^\ast\Biggl(}
{} \times \mathbh{1}\Biggl\{\bigcap_{j=1}^N\bigcap_{u\in A_j} G_K^u \Biggr\}
\mathbh{1}\Biggl\{\bigcap_{j=1}^N\bigcap_{u\in S_j\setminus A_j} [G_K^u]^c
\Biggr\} \Biggr)\\
&&\qquad  \leq\mathop{\sum_{A_i\subseteq S_i}}_{1\leq i\leq N}
\mathbb{E}^\ast\Biggl(\rho^{-\mathcal{N}^{\mathrm{coal}}_\infty\{(\psi
(s),s)\dvtx
s\in\bigcup_{j=1}^N A_j \}}
\mathbh{1}\Biggl\{\bigcap_{j=1}^N\bigcap_{u\in A_j} G_K^u \Biggr\} \Biggr)\nonumber\\
&&\qquad\quad
{} \times\rho^{-\sum_{j=1}^N |S_j\setminus A_j|}
\prod_{j=1}^N\prod_{u\in S_j\setminus A_j}\mathbb{P}^\ast([G_K^u]^c ).\nonumber
\end{eqnarray}
To estimate the expectation in the right-hand side of (\ref{block1}),
we note
that
%
\begin{eqnarray}\label{block2}
\qquad &&\mathcal{N}^{\mathrm{coal}}_\infty\Biggl\{(\psi(s),s)\dvtx s\in\bigcup
_{j=1}^N A_j \Biggr\}\nonumber\\[-8pt]\\[-8pt]
\qquad &&\qquad  \leq\sum_{j=1}^N \mathcal{N}^{\mathrm{coal}}_\infty\{(\psi
(s),s)\dvtx
s\in A_j \}
+\sum_{j=1}^{N-1} \sum_{u\in A_j}
\mathbh{1}\Biggl\{\bigcup_{k=j+1}^N \bigcup_{v\in A_k} M^{u,v} \Biggr\}.\nonumber
\end{eqnarray}
Here we overestimate the number of coalescences of random walks starting
in one ``time-block'' $A_j$ with random walks starting in later
``time-blocks''
$A_k$ by the number of them that meet at least one random walk starting
in a
later ``time-block.'' Together
with H\"older's inequality with $r,r'>1$ and $1/r+1/r'=1$, this
yields
%
\begin{eqnarray}\label{blockholder}
\hspace*{9pt}&&\mathbb{E}^\ast\Biggl( \rho^{-\mathcal{N}^{\mathrm{coal}}_\infty\{(\psi
(s),s)\dvtx
s\in\bigcup_{j=1}^N A_j \}}
\mathbh{1}\Biggl\{\bigcap_{j=1}^N \bigcap_{u \in A_j} G_K^u \Biggr\} \Biggr)\nonumber\\
\hspace*{9pt}&&\qquad \leq\mathbb{E}^\ast\bigl(\rho^{-\sum_{j=1}^N
\mathcal{N}^{\mathrm{coal}}_\infty\{(\psi(s),s)\dvtx s\in A_j \}} \nonumber\\
\hspace*{9pt}&&\qquad \quad\hphantom{\mathbb{E}^\ast\bigl(} {} \times
\rho^{-\sum_{j=1}^{N-1} \sum_{u\in A_j}
\mathbh{1}\{\bigcup_{k=j+1}^N \bigcup_{v\in A_k}
(M^{u,v}\cap G_K^v ) \}} \bigr)\nonumber\\
\hspace*{9pt}&&\qquad \leq\Biggl[\prod_{j=1}^N
\mathbb{E}^\ast\bigl(\rho^{-r \mathcal{N}^{\mathrm{coal}}_\infty\{(\psi
(s),s)\dvtx s\in S_j \}
} \bigr) \Biggr]^{1/r}\\
\hspace*{9pt}&&\qquad \quad {} \times
\Biggl[\mathbb{E}^\ast\Biggl(\prod_{j=1}^{N-1} \prod_{u\in S_j}
\rho^{-r^\prime\mathbh{1}\{\bigcup_{k=j+1}^N\bigcup_{v\in S_k}
M^{u,v}\cap G_K^v \}}
\Biggr) \Biggr]^{1/r'}\nonumber\\
\hspace*{9pt}&&\qquad = \Biggl[\prod_{j=1}^N
\mathbb{E}^\ast\bigl(\rho^{-r \mathcal{N}^{\mathrm{coal}}_\infty\{(\psi
(s),s)\dvtx s\in S_j \}
} \bigr) \Biggr]^{1/r}\nonumber\\
\hspace*{9pt}&&\qquad \quad {} \times
\Biggl[\mathbb{E}^\ast\Biggl(\prod_{j=1}^{N-1} \prod_{u\in S_j}
\Biggl(1+ (\rho^{-r'}-1 ) \mathbh{1}\Biggl\{\bigcup_{k=j+1}^N
\bigcup_{v\in S_k} (M^{u,v}\cap G_K^v ) \Biggr\}
\Biggr) \Biggr) \Biggr]^{1/r'}.\nonumber\hspace*{-30pt}
\end{eqnarray}
In the last step we use the identity
$\rho^{-r^\prime\mathbh{1}\{A\}}=1+(\rho^{-r^\prime}-1)\mathbh
{1}\{A\}$.
Now, by conditional independence and Lemma \ref{Ycondlem}, we have,
for $\epsilon\in(0,(d-4)/2)$ and $1\leq j\leq N-1$,
%
\begin{eqnarray}\label{Ycond}
&&\mathbb{E}^\ast\Biggl( \prod_{u\in S_j} \Biggl(1+ (\rho^{-r'}-1 )
\mathbh{1}\Biggl\{\bigcup_{k=j+1}^N \bigcup_{v\in S_k}
(M^{u,v}\cap G_K^v ) \Biggr\} \Biggr)
\Bigg\vert Y^w, w \in\bigcup_{l>j} S_l \Biggr)\nonumber\\
&&\qquad  \leq\prod_{u\in S_j} \Biggl(1+ (\rho^{-r'}-1 )
\sum_{k=j+1}^N \sum_{v\in S_k}\mathbb{P}^\ast(M^{u,v} \vert
Y^v ) \mathbh{1}\{G_K^v \} \Biggr)\\
&&\qquad  \leq\exp\Biggl[C_\epsilon K (\rho^{-r'}-1 )
\sum_{u\in S_j} \sum_{k=j+1}^N \sum_{v\in S_k}
\frac{1}{(v-u)^{1+\epsilon}} \Biggr].\nonumber
\end{eqnarray}
Clearly,
%
\begin{equation}\label{gapbda}
\sum_{u\in S_j} \sum_{k=j+1}^N \sum_{v\in S_k}\frac
{1}{(v-u)^{1+\epsilon}}
\le\sum_{k=j+1}^N \frac{|S_j| |S_k|}{d(S_j,S_k)^{1+\epsilon}}.
\end{equation}
Substituting this into the right-hand side of (\ref{Ycond}) and using
the resulting
deterministic bounds successively for $j=1,\dots,N-1$, we find that
%
\begin{eqnarray}\label{condbound}
&&\mathbb{E}^\ast\Biggl(\prod_{j=1}^{N-1} \prod_{u\in S_j}
\Biggl(1+ (\rho^{-r'}-1 ) \mathbh{1}\Biggl\{\bigcup_{k=j+1}^N
\bigcup_{v\in S_k} (M^{u,v}\cap G_K^v ) \Biggr\}
\Biggr) \Biggr)\nonumber\\[-8pt]\\[-8pt]
&&\qquad  \leq\exp\biggl[C_\epsilon K (\rho^{-r'}-1 )
\sum_{1\le j < k \le N}
\frac{|S_j| |S_k|}{d(S_j,S_k)^{1+\epsilon}}
\biggr].\nonumber
\end{eqnarray}
It remains to estimate the second factor in the right-hand side of
(\ref{block1}). By
Lemma~\ref{badlm},
%
\begin{equation}\label{secfac}
\rho^{-\sum_{j=1}^N |S_j\setminus A_j|}
\prod_{j=1}^N\prod_{u\in S_j\setminus A_j}\mathbb{P}^\ast([G_K^u]^c )
\le\biggl(\frac{\delta(K)}{\rho} \biggr)^{\sum_{j=1}^N |S_j\setminus A_j|}.
\end{equation}
Observe that, by the binomial formula,
%
\begin{eqnarray}\label{binom}
\mathop{\sum_{A_i\subseteq S_i}}_ {1\leq i\leq N}
\biggl(\frac{\delta(K)}{\rho} \biggr)^{\sum_{j=1}^N |S_j\setminus A_j|}
&=& \biggl(1+\frac{\delta(K)}{\rho} \biggr)^{\sum_{j=1}^N
|S_j|}\nonumber\\[-8pt]\\[-8pt]
&\leq&\exp\Biggl[\frac{\delta(K)}{\rho} \sum_{j=1}^N |S_j|
\Biggr].\nonumber
\end{eqnarray}
Proposition \ref{blockprop} now follows by combining (\ref{block1}) with
(\ref{blockholder}), (\ref{condbound}) and (\ref{secfac}),
and afterward applying (\ref{binom}).
\end{pf*}

\subsection{Continuity at $\kappa=0$}
\label{S4.4}

In this section we prove Theorem \ref{Lyaprop}(ii)(a). We pick $\mu
=\mu_\rho$ as the
starting measure (recall Proposition \ref{lpmunuub}).

By requiring that the $p$ random walks in (\ref{fey-kac2}) do not step
until time $t$, we have, for any $\kappa\in[0,\infty)$,
%
\begin{eqnarray}\label{Lambdalb}
\qquad \Lambda^{\mu_\rho}_p(t;\kappa)
&\geq&\Lambda^{\mu_\rho}_p(t;0)+\frac{1}{pt}
\log\mathrm{P}_0^{\otimes p} \bigl(X_q^\kappa(s)=0 \ \forall
s \in[0,t] \ \forall1\leq q\leq p \bigr)\nonumber\\[-8pt]\\[-8pt]
&=& \Lambda^{\mu_\rho}_p(t;0) - 2d\kappa.\nonumber
\end{eqnarray}
Let $t\to\infty$ to obtain
%
\begin{equation}\label{lambdalb}
\lambda_p(\kappa) \geq\lambda_p(0)- 2d\kappa.
\end{equation}
Therefore, the continuity at $\kappa=0$ reduces to proving that, for all
$d\geq5$, $p\in\mathbb{N}$, $\gamma\in(0,\infty)$ and $\rho\in(0,1)$,
%
\begin{equation}\label{cont0upbd}
\limsup_{\kappa\downarrow0} \lambda_p(\kappa) \leq\lambda_p(0).
\end{equation}

\begin{pf*}{Proof of Theorem \ref{Lyaprop}\textup{(ii)(a)}}
We first give the proof for $p=1$. Fix $L>0$ and $\vartheta\in(0,1)$
arbitrarily. For $j\in\mathbb{N}$, let
%
\begin{eqnarray}\label{cont1}
I_j&=&\bigl[(j-1)L,jL\bigr),\qquad
I_j^\prime=\bigl[(j-1)L,(j-\vartheta)L\bigr),\nonumber\\[-8pt]\\[-8pt]
I_j^{\prime\prime}&=&\bigl[(j-\vartheta)L,jL\bigr)\nonumber
\end{eqnarray}
be the $j$th \textit{time-interval}, \textit{time-block} and \textit{time-gap},
respectively. Fix $r,r^\prime$ with $1/r+1/r^\prime=1$ arbitrarily
and set
%
\begin{equation}\label{mdef}
M = \frac{\rho\gamma(\rho^{-2r^\prime}-1)}{r^\prime\log(1/\rho)}.
\end{equation}
For any Borel set $B\subseteq\mathbb{R}$, let
%
\begin{equation}\label{poissmod}
\widetilde{\Pi}_{\rho\gamma}(B) =
\cases{
\Pi_{\rho\gamma}(B), & \quad if $|\Pi_{\rho\gamma}(B)|\le LM$,\cr
\varnothing, & \quad otherwise.
}
\end{equation}
Since
%
\begin{equation}\label{cont2}
\Pi_{\rho\gamma}([0,t]) \subseteq
\bigcup_{j=1}^{\lceil t/L\rceil} \bigl( \widetilde{\Pi}_{\rho\gamma
}(I_j^\prime)
\cup\bigl(\Pi_{\rho\gamma}(I_j^\prime)
\setminus\widetilde{\Pi}_{\rho\gamma}(I_j^\prime) \bigr)
\cup\Pi_{\rho\gamma}(I_j^{\prime\prime}) \bigr),
\end{equation}
we have
%
\begin{eqnarray}\label{cont3}
\mathcal{N}^{\mathrm{coal}}_\infty\{(X^\kappa(s),s)\dvtx s\in\Pi
_{\rho\gamma
}([0,t]) \}
&\le&\mathcal{N}^{\mathrm{coal}}_\infty\Biggl\{(X^\kappa(s),s)\dvtx
s\in\bigcup_{j=1}^{\lceil t/L\rceil} \widetilde{\Pi}_{\rho\gamma
}(I_j^\prime) \Biggr\}\nonumber\\
&&{} + \sum_{j=1}^{\lceil t/L\rceil} |\Pi_{\rho\gamma}(I_j^\prime)|
\mathbh{1}\{|\Pi_{\rho\gamma}(I_j^\prime)|>LM \}\\
&&{} + \sum_{j=1}^{\lceil t/L\rceil}
|\Pi_{\rho\gamma}(I_j^{\prime\prime})|.\nonumber
\end{eqnarray}
Combining the representation formula (\ref{lyaCRW}) for $p=1$ and
$T=\infty$
with (\ref{cont3}) and applying H\"older's inequality, we find that
%
\begin{equation}\label{cont4}
\exp\bigl[t\bigl(\Lambda^{\mu_\rho}_1(t;\kappa)-\rho\gamma\bigr) \bigr]
\le\mathcal{E}_1 \mathcal{E}_2 \mathcal{E}_3,
\end{equation}
where
\begin{eqnarray}
\label{cont5a}
\qquad \mathcal{E}_1 &=& \bigl( (\mathrm{E}_{ 0} \otimes\mathbb{E}_{\mathrm{Poiss}}\otimes\mathbb{E}^\ast)
\bigl(\rho^{-r\mathcal{N}^{\mathrm{coal}}_\infty\{(X^\kappa(s),s)\dvtx
s\in\bigcup_{j=1}^{\lceil t/L\rceil}
\widetilde{\Pi}_{\rho\gamma}(I_j^\prime) \}}
\bigr) \bigr)^{1/r}, \\
\label{cont5b}
\mathcal{E}_2 &=& \Biggl( \prod_{j=1}^{\lceil t/L\rceil}
\mathbb{E}_{\mathrm{Poiss}}\bigl(\rho^{-r^\prime|\Pi_{\rho\gamma
}(I_j^\prime)|
\mathbh{1}\{|\Pi_{\rho\gamma}(I_j^\prime)|>LM
\}} \bigr) \Biggr)^{1/r^\prime},\\
\label{cont5c}
\mathcal{E}_3 &=& \prod_{j=1}^{\lceil t/L\rceil}
\mathbb{E}_{\mathrm{Poiss}}\bigl( \rho^{-|\Pi_{\rho\gamma}(I_j^{\prime
\prime})|} \bigr)
= \exp[\vartheta(1-\rho)\gamma L\lceil t/L\rceil].
\end{eqnarray}

To estimate $\mathcal{E}_1$ in (\ref{cont5a}), we apply Proposition
\ref{blockprop}
with $\psi(s)=X^\kappa(s)$, $N=\lceil t/L\rceil$,
$S_j=\widetilde{\Pi}_{\rho\gamma}(I_j^\prime)$ and $\rho$
replaced by $\rho^r$.
Then we obtain, for arbitrary $\epsilon\in(0,(d-4)/2)$ and $K>0$,
%
\begin{equation}\label{cont6}
\mathbb{E}^\ast\bigl( \rho^{-r\mathcal{N}^{\mathrm{coal}}_\infty\{
(X^\kappa(s),s)\dvtx
s\in\bigcup_{j=1}^{\lceil t/L\rceil}
\widetilde{\Pi}_{\rho\gamma}(I_j^\prime) \}} \bigr)
\le\mathcal{E}_1^\prime\mathcal{E}_1^{\prime\prime}
\end{equation}
with
%
\begin{equation}\label{cont7}
\mathcal{E}_1^\prime= \Biggl( \prod_{j=1}^{\lceil t/L\rceil} \mathbb
{E}^\ast\bigl(
\rho^{-r^2\mathcal{N}^{\mathrm{coal}}_\infty\{(X^\kappa(s),s)\dvtx
s\in\widetilde{\Pi}_{\rho\gamma}(I_j^\prime) \}}
\bigr) \Biggr)^{1/r}
\end{equation}
and
%
\begin{eqnarray}\label{cont8}
\mathcal{E}_1^{\prime\prime}
&=& \exp\Biggl[ \frac{\delta(K)}{\rho^r}
\sum_{j=1}^{\lceil t/L\rceil}
|\widetilde{\Pi}_{\rho\gamma}(I_j^\prime)|\nonumber\\[-8pt]\\[-8pt]
&&\hphantom{\exp\Biggl[}{}
+ C_\epsilon K\frac{\rho^{-rr^\prime}-1}{r^\prime}
\sum_{1\le j<k\le\lceil t/L\rceil}
\frac{|\widetilde{\Pi}_{\rho\gamma}(I_j^\prime)|
|\widetilde{\Pi}_{\rho\gamma}(I_k^\prime)|
}{d(I_j,I_k)^{1+\epsilon}}\Biggr].\nonumber
\end{eqnarray}
To estimate $\mathcal{E}_1^\prime$, we write
%
\begin{equation}\label{cont9}
\Pi_{\rho\gamma} = \Pi_{\rho^{r^2}\gamma}^{(1)}
\cup\Pi_{(\rho-\rho^{r^2})\gamma}^{(2)},
\end{equation}
where $\Pi_{\rho^{r^2}\gamma}^{(1)}$ and $\Pi_{(\rho-\rho
^{r^2})\gamma}^{(2)}$
are independent Poisson processes on $\mathbb{R}$ with intensity $\rho
^{r^2}\gamma$
and $(\rho-\rho^{r^2})\gamma$, respectively, and we use that
[recall (\ref{poissmod})]
%
\begin{eqnarray}\label{cont10}
&&\mathcal{N}^{\mathrm{coal}}_\infty\{(X^\kappa(s),s)\dvtx
s\in\widetilde{\Pi}_{\rho\gamma}(I_j^\prime)
\}\nonumber\\[-8pt]\\[-8pt]
&&\qquad  \le\mathcal{N}^{\mathrm{coal}}_\infty\bigl\{(X^\kappa(s),s)\dvtx
s\in\Pi_{\rho^{r^2}\gamma}^{(1)}(I_j) \bigr\}
+ \bigl|\Pi_{(\rho-\rho^{r^2})\gamma}^{(2)}(I_j)\bigr|.\nonumber
\end{eqnarray}
This leads to
%
\begin{eqnarray}\label{cont11}
\mathcal{E}_1^\prime&\le& \Biggl( \prod_{j=1}^{\lceil t/L\rceil}
\mathbb{E}^\ast\bigl( \rho^{-r^2\mathcal{N}^{\mathrm{coal}}_\infty\{
(X^\kappa(s),s)\dvtx
s\in\Pi_{\rho^{r^2}\gamma}(I_j) \}}
\bigr) \Biggr)^{1/r} \nonumber\\[-8pt]\\[-8pt]
&&{}\times\exp\biggl[ (\rho-\rho^{r^2})\gamma\frac{\rho^{-r^2}-1}{r}
L\lceil t/L\rceil\biggr].\nonumber
\end{eqnarray}
To estimate $\mathcal{E}_1^{\prime\prime}$, note that
$|\widetilde{\Pi}_{\rho\gamma}(I_j^\prime)|\le LM$ for all $j$
and $d(I_j^\prime,I_k^\prime)\ge\vartheta L(k-j)$ for $k>j$, so that
%
\begin{equation}\label{cont12}
\mathcal{E}_1^{\prime\prime} \le\exp\biggl[ \biggl(\frac{\delta(K)}{\rho^r}M
+ C_\epsilon^\prime K\frac{\rho^{-rr^\prime}-1}{r^\prime}
\frac{M^2}{\vartheta^{1+\epsilon}L^\epsilon} \biggr) L\lceil t/L\rceil
\biggr],
\end{equation}
where $C_\epsilon^\prime=C_\epsilon\sum_{j=1}^\infty
j^{-(1+\epsilon)}$.
Since the distribution of $\mathcal{N}^{\mathrm{coal}}_\infty$ is
invariant w.r.t. spatial
shifts of the coalescing random walks, and $X^\kappa$ and $\Pi_{\rho
^{r^2}\gamma}$
have independent and stationary increments, we obtain
%
\begin{eqnarray}\label{cont13}
&&(\mathrm{E}_{ 0}\otimes\mathbb{E}_{\mathrm{Poiss}}) \Biggl( \prod
_{j=1}^{\lceil t/L\rceil}
\mathbb{E}^\ast\bigl( \rho^{-r^2\mathcal{N}^{\mathrm{coal}}_\infty\{
(X^\kappa(s),s)\dvtx
s\in\Pi_{\rho^{r^2}\gamma}(I_j) \}} \bigr) \Biggr) \nonumber\\
&&\qquad  = (\mathrm{E}_{ 0} \otimes\mathbb{E}_{\mathrm{Poiss}}) \Biggl( \prod
_{j=1}^{\lceil t/L\rceil}
\mathbb{E}^\ast\bigl( \rho^{-r^2\mathcal{N}^{\mathrm{coal}}_\infty\{
(X^\kappa(s)-X^\kappa
((j-1)L),s)\dvtx
s\in\Pi_{\rho^{r^2}\gamma}(I_j) \}} \bigr)
\Biggr)\nonumber\\[-8pt]\\[-8pt]
&&\qquad  = \bigl( (\mathrm{E}_{ 0} \otimes\mathbb{E}_{\mathrm{Poiss}}\otimes
\mathbb{E}^\ast)
\bigl( \rho^{-r^2\mathcal{N}^{\mathrm{coal}}_\infty\{(X^\kappa
(s),s)\dvtx
s\in\Pi_{\rho^{r^2}\gamma}([0,L]) \}} \bigr)
\bigr)^{\lceil t/L\rceil}\nonumber \\
&&\qquad  = \exp\bigl[ \bigl(\Lambda_1^{\mu_{\rho^{r^2}}}(L;\kappa)-\rho
^{r^2}\gamma\bigr)
L\lceil t/L\rceil\bigr],\nonumber
\end{eqnarray}
where in the last line we have used the representation formula (\ref{lyaCRW})
for $p=1$, $T=\infty$ and $\rho$ and $t$ replaced by $\rho^{r^2}$
and $L$,
respectively. Now substitute (\ref{cont11}) and (\ref{cont12}) into
(\ref{cont6}), substitute the obtained inequality into (\ref{cont5a}) and
use (\ref{cont13}) to arrive at
%
\begin{eqnarray}\label{cont14}
\mathcal{E}_1 &\le&\exp\biggl[ \frac{1}{r^2} \bigl(\Lambda_1^{\mu_{\rho
^{r^2}}}(L;\kappa)
-\rho^{r^2}\gamma\bigr) L\lceil t/L\rceil\biggr] \nonumber\\
&&{} \times\exp\biggl[ \biggl( (\rho-\rho^{r^2})\gamma\frac{\rho^{-r^2}-1}{r^2}
+ \frac{\delta(K)}{r\rho^r}M \\
&&\hspace*{66pt}
{} + C_\epsilon^\prime K\frac{\rho^{-rr^\prime}-1}{rr^\prime}
\frac{M^2}{\vartheta^{1+\epsilon}L^\epsilon} \biggr) L\lceil t/L\rceil\biggr].\nonumber
\end{eqnarray}

We next estimate $\mathcal{E}_2$ in (\ref{cont5b}). Using Chebyshev's
exponential
inequality, we obtain, for $j=1,\ldots,\lceil t/L\rceil$,
%
\begin{eqnarray}\label{cont15}
&&\mathbb{E}_{\mathrm{Poiss}}\bigl( \rho^{-r^\prime|\Pi_{\rho\gamma
}(I_j^\prime)|
\mathbh{1}\{|\Pi_{\rho\gamma}(I_j^\prime)|>LM \}}
\bigr) \nonumber\\
&&\qquad  \le1 + \mathbb{E}_{\mathrm{Poiss}}\bigl( \rho^{-r^\prime|\Pi_{\rho
\gamma
}(I_j^\prime)|}
\mathbh{1}\{|\Pi_{\rho\gamma}(I_j^\prime)|>LM \} \bigr) \nonumber\\
&&\qquad  \le1 + \rho^{r^\prime LM}
\mathbb{E}_{\mathrm{Poiss}}\bigl( \rho^{-2r^\prime|\Pi_{\rho\gamma
}(I_j^\prime)|} \bigr) \\
&&\qquad  \le1 + \rho^{r^\prime LM}
\mathbb{E}_{\mathrm{Poiss}}\bigl( \rho^{-2r^\prime|\Pi_{\rho\gamma
}(I_j)|} \bigr) \nonumber\\
&&\qquad  = 1 + \exp\bigl[ \bigl(\rho\gamma(\rho^{-2r^\prime}-1)
- r^\prime M\log(1/\rho) \bigr) L \bigr].\nonumber
\end{eqnarray}
By our choice of $M$ in (\ref{mdef}), the expression in the right-hand side
equals $2$, and we conclude that
%
\begin{equation}\label{cont16}
\mathcal{E}_2 \le e^{\lceil t/L\rceil}.
\end{equation}

Finally, substitute (\ref{cont14}), (\ref{cont16}) and (\ref
{cont5c}) into
(\ref{cont4}), take the logarithm on both sides of the resulting inequality,
divide by $t$, pass to the limit as $t\to\infty$ and recall (\ref{lyap2}).
Then we obtain
%
\begin{eqnarray}\label{cont17}
\qquad\quad  \lambda_1^{\mu_\rho}(\kappa)-\rho\gamma&\le&\frac{1}{r^2}
\bigl( \Lambda_1^{\mu_{\rho^{r^2}}}(L;\kappa)-\rho^{r^2}\gamma\bigr)
+ (\rho-\rho^{r^2})\gamma\frac{\rho^{-r^2}-1}{r^2}
\nonumber\\[-8pt]\\[-8pt]
&&{} + \frac{\delta(K)}{r\rho^r}M
+ C_\epsilon^\prime K\frac{\rho^{-rr^\prime}-1}{rr^\prime}
\frac{M^2}{\vartheta^{1+\epsilon}L^\epsilon} + \frac{1}{L}
+ \vartheta(1-\rho)\gamma.\nonumber
\end{eqnarray}
As can be seen from (\ref{fey-kac2}),
$\kappa\mapsto\Lambda_1^{\mu_{\rho^{r^2}}}(L;\kappa)$ is
continuous at
$\kappa=0$. Hence, passing in (\ref{cont17}) to the limits as
$\kappa\downarrow0$, $L\to\infty$, $K\to\infty$ and $\vartheta
\downarrow0$
(in this order), we find that
%
\begin{equation}\label{cont18}
\qquad \limsup_{\kappa\downarrow0} \bigl(\lambda_1^{\mu_\rho}(\kappa)-\rho
\gamma\bigr)
\le\frac{1}{r^2} \bigl(\lambda_1^{\mu_{\rho^{r^2}}}(0)-\rho
^{r^2}\gamma\bigr)
+ (\rho-\rho^{r^2})\gamma\frac{\rho^{-r^2}-1}{r^2}.
\end{equation}
Expanding the exponential function in the right-hand side of (\ref{fey-kac2})
into a Taylor series and using (\ref{dual-vm}), we see that
$\rho\mapsto\Lambda_1^{\mu_\rho}(t;0)$ is nondecreasing. Hence,
the same is
true for $\rho\mapsto\lambda_1^{\mu_\rho}(0)$. Taking this into account,
we may finally pass to the limit as $r\downarrow1$ in (\ref{cont18})
to arrive at
%
\begin{equation}\label{cont19}
\limsup_{\kappa\downarrow0} \bigl(\lambda_1^{\mu_\rho}(\kappa)-\rho
\gamma\bigr)
\le\lambda_1^{\mu_\rho}(0)-\rho\gamma.
\end{equation}
This is the desired inequality (\ref{cont0upbd}) for $p=1$.

The extension to $p\in\mathbb{N}\setminus\{1\}$ is straightforward. The
proof follows
the same arguments with $X^\kappa$ and $\Pi_{\rho\gamma}$ replaced
by $p$
independent copies $X^\kappa_q$ and $\Pi_{\rho\gamma}^{(q)}$,
$q=1,\ldots,p$, of $X^\kappa$ and $\Pi_{\rho\gamma}$, respectively.
\end{pf*}

\subsection{Large $\kappa$}
\label{S4.5}

In this section we prove Theorem \ref{Lyaprop}(ii)(b). We again pick
$\mu=\mu_\rho$ as the starting measure (recall Proposition \ref{lpmunuub}).

\begin{pf*}{Proof of Theorem \ref{Lyaprop}\textup{(ii)(b)}}
Recall (\ref{lyasand}). We first give the proof for \mbox{$p=1$}. We show
that, for all $\rho\in(0,1)$, $\gamma>0$ and $L>0$,
%
\begin{equation}\label{lkappa1}
\lim_{\kappa\to\infty} \Lambda_1^{\mu_\rho}(L;\kappa) = \rho
\gamma.
\end{equation}
Then the claim for $p=1$ follows from (\ref{cont17}) by passing to the limits
as $\kappa\to\infty$, $L\to\infty$, $K\to\infty$, $\vartheta
\downarrow0$ and
$r\downarrow1$ (in this order).

To prove (\ref{lkappa1}), we use the representation formula (\ref{lyaCRW}):
%
\begin{eqnarray}\label{lkappa2}
&&\Lambda_1^{\mu_\rho}(L;\kappa)-\rho\gamma\nonumber\\[-8pt]\\[-8pt]
&&\qquad = \frac{1}{L} \log
(\mathrm{E}_{ 0} \otimes\mathbb{E}_{\mathrm{Poiss}}\otimes\mathbb
{E}^\ast)
\bigl( \rho^{-\mathcal{N}^{\mathrm{coal}}_\infty\{(X^\kappa(s),s)\dvtx
s\in\Pi_{\rho\gamma}([0,L]) \}} \bigr).\nonumber
\end{eqnarray}
Recall that we are in a transient situation ($d\ge5$) and write
$X^\kappa(s)=X^1(\kappa s)$. Then, $\mathrm{P}_{ 0}\otimes\mathbb
{P}_{\mathrm{Poiss}}$-a.s.
%
\begin{equation}\label{lkappa3}
\lim_{\kappa\to\infty} \mathop{\min_{s_1,s_2\in\Pi_{\rho\gamma
}([0,L])}}_{ s_1\ne s_2}
|X^\kappa(s_1)-X^\kappa(s_2) | = \infty,
\end{equation}
and, consequently,
%
\begin{equation}\label{lkappa4}
\qquad\hspace*{7pt} \lim_{\kappa\to\infty} \mathcal{N}^{\mathrm{coal}}_\infty\{
(X^\kappa(s),s)\dvtx
s\in\Pi_{\rho\gamma}([0,L]) \} = 0\qquad
\mbox{in probability w.r.t. }\mathbb{P}^\ast.\hspace*{-7pt}
\end{equation}
Since, moreover, $\mathcal{N}^{\mathrm{coal}}_\infty\{(X^\kappa
(s),s)\dvtx
s\in\Pi_{\rho\gamma}([0,L]) \} \le|\Pi_{\rho\gamma}([0,L])|$,
we may apply Lebesgue's dominated convergence theorem to see that
the expression on the right of (\ref{lkappa2}) converges to $0$ as
$\kappa\to\infty$. This proves (\ref{lkappa1}).

The extension to $p\in\mathbb{N}\setminus\{1\}$ is easy. Indeed,
by (\ref{fey-kac1})--(\ref{lyapdef}) and Jensen's inequality,
%
\begin{eqnarray}\label{Lamprel}
\exp[pt\Lambda^{\mu_\rho}_p(t;\kappa,\gamma) ]
&=&\mathbb{E}_{ \mu_\rho} \biggl( \biggl[\mathrm{E}_{ 0} \biggl(\exp\biggl[
\gamma\int_0^t\xi\bigl(X^\kappa(s),t-s \bigr) \,ds \biggr] \biggr) \biggr]^p \biggr)\nonumber\\
&\leq&\mathbb{E}_{ \mu_\rho} \biggl(\mathrm{E}_{ 0} \biggl(\exp\biggl[
p\gamma\int_0^t\xi\bigl(X^\kappa(s),t-s \bigr) \,ds \biggr] \biggr) \biggr)\\
&=& \exp[t\Lambda^{\mu_\rho}_1(t;\kappa,p\gamma) ].\nonumber
\end{eqnarray}
Let $t\to\infty$ to get
%
\begin{equation}\label{lamprel}
\lambda_p(\kappa;\gamma) \leq\frac{1}{p} \lambda_1(\kappa
;p\gamma).
\end{equation}
This together with the assertion for $p=1$ and (\ref{lyasand})
implies the claim for arbitrary $p\in\mathbb{N}$.
\end{pf*}


\section[Proof of Theorem 1.5(i) and (ii)(c)]{Proof of Theorem \protect\ref{Lyaprop}(i) and (ii)(c)}
\label{S5}

Throughout this section we assume that $p(\cdot,\cdot)$ satisfies
(\ref{pdef}) and has
zero mean and finite variance. Theorem \ref{Lyaprop}(i) is proved in
Section \ref{S5.1}
and Theorem \ref{Lyaprop}(ii)(c) in Section \ref{S5.2}. As a starting measure
we pick $\mu=\nu_\rho$ (recall Proposition \ref{lpmunuub}).


\subsection{Triviality in low dimensions}
\label{S5.1}

The proof of Theorem \ref{Lyaprop}(i) is similar to that of Theorem 1.3.2(i)
in G\"artner, den Hollander and Maillard \cite{garholmai07}. The key
observation is the following:

\begin{lemma}\label{modld}
If $1\leq d\leq4$, then for any finite $Q\subset\mathbb{Z}^d$ and
$\rho\in(0,1)$,
%
\begin{equation}\label{stirineq3}
\lim_{t\to\infty} \frac{1}{t} \log\mathbb{P}_{\nu_\rho} \bigl(\xi
(x,s) = 1
\ \forall x\in Q \ \forall s\in[0,t] \bigr) = 0.
\end{equation}
\end{lemma}

\begin{pf}
In the spirit of Bramson, Cox and Griffeath \cite{bramcoxgri88},
Section 1,
we argue as follows. The graphical representation of the VM (recall Section
\ref{S1.3}) allows us to write down a suitable expression for the probability
in (\ref{stirineq3}). Indeed, let
%
\begin{equation}\label{stirdef}
H_t^Q = \{x\in\mathbb{Z}^d\dvtx\mbox{ there is a path from } (x,0)
\mbox{ to } Q \times[0,t] \mbox{ in } {\mathcal{G}}_t \},
\end{equation}
where, as in Section \ref{S1.3}, $\mathcal{G}_t$ is the graphical
representation of the
voter model up to time $t$ (see Figure \ref{fig4}).


%
%
%
%
%
%
%
%
%
%

\begin{figure}[b]

\includegraphics{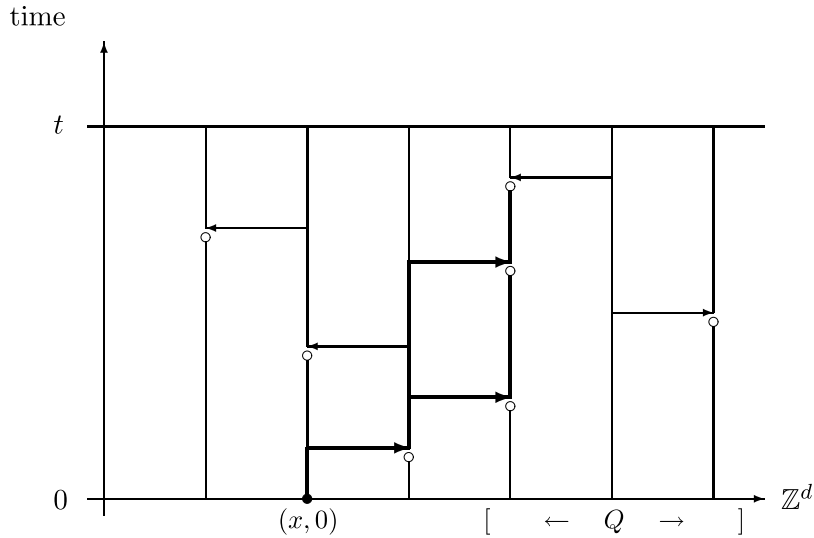}

  \caption{Some paths from $(x,0)$ to $Q\times[0,t]$ in ${\mathcal{G}}_t$.}\label{fig4}
\end{figure}


Note that $H^Q_0=Q$ and that $t\mapsto H^Q_t$ is nondecreasing.
Denote by $\mathcal{P}$ and $\mathcal{E}$, respectively,
probability and expectation associated with the graphical representation
$\mathcal{G}_t$. Then
%
\begin{equation}\label{stirrep}
\mathbb{P}_{\nu_\rho} \bigl(\xi(x,s) = 1
\ \forall x\in Q \ \forall s\in[0,t] \bigr)
= ({\mathcal{P}}\otimes\nu_\rho) \bigl(H_t^Q \subseteq\xi(0) \bigr),
\end{equation}
where $\xi(0)=\{x\in\mathbb{Z}^d\dvtx\xi(x,0)=1\}$ is the set of initial
locations of
$1$'s. Indeed, (\ref{stirrep}) holds because if $\xi(x,0)=0$ for some
$x \in H^Q_t$,
then this 0 will propagate into $Q$ prior to time $t$ (see Figure \ref{fig4}).

By Jensen's inequality,
%
\begin{equation}\label{stirineq1}
({\mathcal{P}}\otimes\nu_\rho) \bigl(H_t^Q \subseteq\xi(0) \bigr)
= {\mathcal{E}} \bigl(\rho^{|H_t^Q|} \bigr)
\geq\rho^{ {\mathcal{E}}|H_t^Q|}.
\end{equation}
Moreover, $H_t^Q = \bigcup_{y\in Q} H_t^{\{y\}}$, implying
%
\begin{equation}\label{stirineq2}
{\mathcal{E}}|H_t^Q| \leq|Q| {\mathcal{E}}\bigl|H_t^{\{0\}}\bigr|.
\end{equation}
By the dual graphical representation, $|H_t^{\{0\}}|$ coincides
in distribution with the number of coalescing random walks alive at
time $t$
when starting at site $0$ at times generated by a rate $1$ Poisson stream.
As shown in Bramson, Cox and Griffeath \cite{bramcoxgri88}, Theorem 2,
if $p(\cdot,\cdot)$ is a simple random walk, then
%
\begin{equation}\label{limE}
{\mathcal{E}}\bigl|H_t^{\{0\}}\bigr|=o(t)
\qquad \mbox{as } t\to\infty\mbox{ when } 1\leq d\leq4,
\end{equation}
in which case (\ref{stirineq3}) follows from (\ref{stirrep})--(\ref
{stirineq2}).
As noted in Bramson, Cox and Le Gall \cite{bramcoxleg01}, Lemma 2, and
its proof, the key
ingredient in the proof of (\ref{limE}) extends from a simple random
walk to
a random walk with zero mean and finite variance.
\end{pf}

We are now ready to give the proof of Theorem \ref{Lyaprop}(i).

\begin{pf*}{Proof of Theorem \ref{Lyaprop}\textup{(i)}}
Fix $1\leq d \leq4$, $\kappa\in[0,\infty)$, $\gamma\in(0,\infty
)$ and $\rho\in(0,1)$.
Since $p\mapsto\lambda_p(\kappa)$ is nondecreasing and $\lambda
_p(\kappa)\leq\gamma$
for all $p\in\mathbb{N}$ [recall (\ref{lyasand})], it suffices to
give the
proof for $p=1$.
For $p=1$, (\ref{fey-kac2}) reads
%
\begin{equation}\label{L1a}
\Lambda^{\nu_\rho}_1(t) = \frac{1}{t} \log(\mathbb{E}_{ \nu_\rho}
\otimes\mathrm{E}_{ 0})
\biggl(\exp\biggl[\gamma\int_0^t
\xi\bigl(X^\kappa(s),t-s \bigr) \,ds \biggr] \biggr).
\end{equation}
By restricting $X^\kappa$ to stay inside a finite box $Q\subset
\mathbb{Z}^d$
around 0 up
to time $t$ and requiring $\xi$ to be 1 in the entire box up to time
$t$, we
obtain
%
\begin{eqnarray}\label{L1b}
&&(\mathbb{E}_{ \nu_\rho} \otimes\mathrm{E}_{ 0}) \biggl(\exp\biggl[\gamma
\int_0^t
\xi\bigl(X^\kappa(s),t-s\bigr) \,ds \biggr]
\biggr)\nonumber\\[-8pt]\\[-8pt]
&&\qquad  \geq e^{\gamma t} \mathbb{P}_{\nu_\rho} \bigl(\xi(x,s) = 1
\ \forall x\in Q \ \forall s \in[0,t] \bigr)
\mathrm{P}_0 \bigl(X^\kappa(s)\in Q \ \forall s\in[0,t] \bigr).\nonumber
\end{eqnarray}
The first factor is $e^{o(t)}$ by Lemma \ref{modld}. For the second
factor, we have
%
\begin{equation}\label{L1d}
\lim_{t\to\infty} \frac{1}{t} \log
\mathrm{P}_0 \bigl(X^\kappa(s)\in Q \ \forall s\in[0,t] \bigr)
= \lambda^\kappa(Q),
\end{equation}
with $\lambda^\kappa(Q)<0$ the principal Dirichlet eigenvalue on $Q$
of $\kappa\Delta$,
the generator of $X^\kappa$. Combining (\ref{stirineq3}) and (\ref
{L1a})--(\ref{L1d}), we
arrive at
%
\begin{equation}\label{resSE}
\lambda_1(\kappa) = \lim_{t\to\infty} \Lambda^{\nu_\rho}_1(t)
\geq\gamma+\lambda^\kappa(Q).
\end{equation}
Finally, let $Q\uparrow\mathbb{Z}^d$ and use that $\lim_{Q\uparrow
\mathbb{Z}
^d}\lambda^\kappa(Q)=0$
(see, e.g., Spitzer \cite{sp76}, Section~21) to arrive
at $\lambda_1(\kappa)\geq\gamma$. Since, trivially, $\lambda
_1(\kappa)\leq\gamma$,
we get $\lambda_1(\kappa)=\gamma$.
\end{pf*}


\subsection{Intermittency for small $\kappa$}
\label{S5.2}

We start this section by recalling some large deviation results for the
VM that will be
needed to prove Theorem \ref{Lyaprop}(ii)(c). Cox and Griffeath \cite
{coxgri83} showed
that for the VM with a simple random walk transition kernel given by
(\ref{SRW}), the
occupation time of the origin up to time $t \geq0$,
%
\begin{equation}\label{occtime}
T_t=\int_0^t \xi(0,s)\, ds,
\end{equation}
satisfies a strong law of large numbers and a central limit theorem for
$d\geq2$. For $d=1$ there is no law of large numbers: $T_t/t$ has a
nontrivial limiting law. These results carry over to a random walk with
zero mean and
finite variance.

The following proposition gives \textit{large deviation bounds}.

\begin{proposition}[(Bramson, Cox and Griffeath \cite{bramcoxgri88}, Theorem 1;
Bramson, Cox
and Le Gall \cite{bramcoxleg01}, Lemma 2 and its proof; Maillard and
Mountford \cite{mailmoun07}, Theorem~1.3.2)]\label{LDbds}
Suppose that $p(\cdot,\cdot)$ has zero mean and finite variance. Then
for every
$\alpha\in(\rho,1)$ there exist $0<I^-(\alpha)<I^+(\alpha)<\infty
$ such that,
for $t$ sufficiently large (depending on $\alpha$),
%
\begin{equation}\label{ratevm}
e^{-I^+(\alpha) b_t} \leq\mathbb{P}_{\nu_\rho} \biggl(\frac{T_t}t\geq
\alpha\biggr)
\leq e^{-I^-(\alpha) b_t}
\end{equation}
with
%
\begin{equation}\label{ratevm2}
b_t =
\cases{
\log t,
&\quad if $ d=2$,\cr
\sqrt{t},
&\quad if  $d=3$,\cr
\displaystyle\frac{t}{\log t},
&\quad if  $d=4$,\cr
t,
&\quad if  $d \geq5$.
}
\end{equation}
\end{proposition}

By interchanging the opinions $0$ and $1$, similar bounds are obtained
for $\mathbb{P}_{\nu_\rho}
(T_t/t\leq\alpha)$, $\alpha\in(0,\rho)$. The case $\alpha=1$ may
be included in
$d \geq3$ but not in $d=2$, for which it is shown in Maillard and
Mountford \cite{mailmoun07},
Theorem 1.3.1, that $\mathbb{P}(T_t=t)$ is of order $\exp[-(\log
t)^2]$. A
full large deviation principle
is expected to hold for $d\geq3$, but this has not been established.
Inspection of the proof
in Bramson, Cox and Griffeath \cite{bramcoxgri88} shows that for $d
\geq5$ there exists a
$C>0$ such that
%
\begin{equation}\label{Iminlb}
I^-(\alpha)\geq C \bigl(\sqrt{\alpha}-\sqrt{\rho} \bigr)^2,\qquad  \alpha\in(\rho,1).
\end{equation}
No comparable upper bound on $I^+$ is given.


We are now ready to give the proof of Theorem \ref{Lyaprop}(ii)(c).

\begin{pf*}{Proof of Theorem \ref{Lyaprop}\textup{(ii)(c)}}
We first give the proof for $\kappa=0$. Fix $d \geq5$, $p\in\mathbb{N}$,
$\gamma\in(0,\infty)$ and $\rho\in(0,1)$, and recall that
$\lambda_p(0)>\rho\gamma$ by Theorem \ref{Lyalip}(ii). Pick
$\alpha\in(\rho,\gamma^{-1}\lambda_p(0))$ and define
%
\begin{equation}\label{int-3}
I(\alpha)=-\limsup_{t\to\infty}\frac1t \log\mathbb{P}_{ \nu
_\rho}
\biggl(\frac{1}{t}T_t\ge\alpha\biggr) > 0,
\end{equation}
where the positivity of the limit comes from the upper bound in (\ref{ratevm}),
which implies $I(\alpha) \geq I^-(\alpha)>0$. Put
%
\begin{equation}\label{int-5}
\beta=\gamma^{-1} \biggl[\lambda_p(0)+\frac{1}{2p}I(\alpha) \biggr]
\end{equation}
and split
%
\begin{equation}\label{int-7}
\Lambda^{\nu_\rho}_p(t)
=\frac1{pt}\log\mathbb{E}_{ \nu_\rho} (e^{ p\gamma T_t} )
=\frac1{pt}\log(A_t+B_t+C_t)
\end{equation}
with
%
\begin{eqnarray}\label{int-9}
A_t &= &\mathbb{E}_{ \nu_\rho} \biggl(e^{ p\gamma T_t} \mathbh{1}\biggl\{0\leq
\frac1tT_t<\alpha\biggr\} \biggr),\nonumber\\
B_t &=& \mathbb{E}_{ \nu_\rho} \biggl(e^{ p\gamma T_t} \mathbh{1}\biggl\{\alpha
\leq\frac1tT_t<\beta\biggr\} \biggr),\\
C_t &=& \mathbb{E}_{ \nu_\rho} \biggl(e^{ p\gamma T_t} \mathbh{1}\biggl\{\frac
1tT_t\geq\beta\biggr\} \biggr).\nonumber
\end{eqnarray}
Next, note that
%
\begin{equation}\label{int-13}
A_t \leq e^{ p\gamma\alpha t},\qquad  B_t \leq e^{ p\gamma\beta t}
\mathbb{P}_{ \nu_\rho} \biggl(\frac1tT_t\geq\alpha\biggr).
\end{equation}
Thus, in (\ref{int-7}) both $A_t$ and $B_t$ are negligible as $t\to
\infty$, because
$\lim_{t\to\infty}\Lambda^{\nu_\rho}_p(t)=\lambda_p(0)$ while
$\gamma\alpha<\lambda_p(0)$
and $\gamma\beta-\frac1p I(\alpha)=\lambda_p(0)-\frac
{1}{2p}I(\alpha)<\lambda_p(0)$. Hence,
%
\begin{equation}\label{int-17}
\lambda_p(0)=\lim_{t\to\infty}\frac1{pt}\log C_t.
\end{equation}
Now, by (\ref{int-5}) and (\ref{int-17}), we have
%
\begin{eqnarray}\label{int-20}
\lambda_{p+1}(0)
&=&\lim_{t\to\infty}\frac1{(p+1)t} \log\mathbb{E}_{ \nu_\rho}
\bigl(e^{(p+1)\gamma T_t} \bigr)\nonumber\\
&\geq&\limsup_{t\to\infty}\frac1{(p+1)t} \log\mathbb{E}_{ \nu
_\rho}
\biggl(e^{(p+1)\gamma T_t} \mathbh{1}\biggl\{\frac1tT_t\geq\beta\biggr\}
\biggr)\nonumber\\
&\geq&\frac1{p+1} \gamma\beta+\lim_{t\to\infty}\frac1{(p+1)t}
\log C_t\\
&=&\frac1{p+1} \gamma\beta+\frac{p}{p+1} \lambda_p(0)\nonumber\\
&=&\lambda_p(0)+\frac{1}{2p(p+1)} I(\alpha)>\lambda_p(0),\nonumber
\end{eqnarray}
which proves the gap between $\lambda_p(0)$ and $\lambda_{p+1}(0)$.

By the continuity of $\kappa\mapsto\lambda_p(\kappa)$ at $\kappa
=0$ in
Theorem \ref{Lyaprop}(ii)(a), it follows that there exists $\kappa_0>0$
such that $\lambda_p(\kappa)>\lambda_{p-1}(\kappa)$ for all $\kappa
\in[0,\kappa_0)$
when $p=2$ and, by the remark made after formula (\ref{lyasand}), also when
$p\in\mathbb{N}\setminus\{1\}$.
\end{pf*}


\begin{appendix}

\section*{\texorpdfstring{Appendix: Heuristic explanation of
Conjecture \protect\lowercase{\ref{cj-larg-kappa}}}{Appendix: Heuristic explanation of Conjecture 1.6}}
\label{appA}
\setcounter{thm}{0}
\setcounter{equation}{0}

In this appendix we give a heuristic explanation of (\ref{limlamb*}).
We only
consider the case $p=1$. A similar argument works for $p\in\mathbb{N}
\setminus\{1\}$.
As starting measure we pick $\mu=\mu_\rho$ (recall Proposition \ref{lpmunuub}).

1. \textit{Pair correlation}.
Lemma \ref{reprform} for $n=2$ yields the following representation
for the pair correlation function of the VM in equilibrium.

\begin{lemma}\label{2corr}
Suppose that $p(\cdot,\cdot)$ is symmetric and transient. Then, for
all $x_1,x_2\in\mathbb{Z}^d$
and $s\geq0$,
%
\begin{equation}\label{2corr-1}
\mathbb{E}_{ \mu_\rho} \bigl([\xi(x_1,s)-\rho][\xi(x_2,0)-\rho] \bigr)
=\frac{\rho(1-\rho)}{G_d} \int_0^\infty p_{s+t}(x_1,x_2) \,dt
\end{equation}
with $G_d = \int_0^\infty p_t(0,0) \,dt$.
\end{lemma}

\begin{pf}
The proof is standard. By (\ref{dual-vm}) with $T=\infty$ and $n=2$,
we have
%
\begin{eqnarray}\label{AP-1}
&&\mathbb{E}_{ \mu_\rho} \bigl([\xi(x_1,s)-\rho][\xi(x_2,0)-\rho]
\bigr)\nonumber\\[-8pt]\\[-8pt]
&&\qquad = \rho(1-\rho) \mathbb{P}^\ast\bigl(\mathcal{N}_\infty\{
(x_1,0),(x_2,s)\}=1 \bigr).\nonumber
\end{eqnarray}
The probability in the right-hand side of (\ref{AP-1}) can be computed
as follows.
The first random walk starts from site $x_1$ at time $0$, moves freely
until time
$s$, and reaches some site $y$ at time $s$. The second random walk
starts from
site $x_2$ at time $s$ and has to eventually coalesce with the first
random walk.
This gives
%
\begin{equation}\label{coal2-5*}
\mathbb{P}^\ast\bigl(\mathcal{N}_\infty\{(x_1,0),(x_2,s)\}=1 \bigr)
= \sum_{y\in\mathbb{Z}^d} p_s(x_1,y) w(y-x_2)
\end{equation}
with
%
\begin{equation}\label{coal2-7}
w(z) = \mathbb{P}_z (Z_t=0 \mbox{ for some } 0 \leq t < \infty),
\qquad z\in\mathbb{Z}^d.
\end{equation}
Here we use that, by the symmetry of $p(\cdot,\cdot)$, the difference between
the two random walks is a single random walk $Z$ running at double the
speed. By a
renewal argument (see Spitzer \cite{sp76}, Section 4), for transient
$p(\cdot,\cdot)$
we have
%
\begin{equation}\label{coal2-ex}
w(z) = \frac{1}{G_d} \int_0^\infty p_t(z,0) \,dt.
\end{equation}
Combining (\ref{AP-1}), (\ref{coal2-5*}) and (\ref{coal2-ex}), we
obtain (\ref{2corr-1}).
\end{pf}

2. \textit{Green term}.
From now on let $p(\cdot,\cdot)$ be a simple random walk. Fix $d \geq
5$, $\gamma\in(0,\infty)$
and $\rho\in(0,1)$. Scaling time by $\kappa$ in (\ref{fey-kac2}),
we have
$\lambda_1(\kappa)=\kappa\lambda_1^*(\kappa)$ with
%
\begin{equation}\label{lambda*scal*}
\lambda_1^*(\kappa) = \lim_{t\to\infty} \Lambda_1^*(\kappa;t)
\end{equation}
and
%
\begin{equation}\label{lambda*scal}
\Lambda_1^*(\kappa;t)
= \frac{1}{t} \log(\mathbb{E}_{ \mu_\rho}\otimes\mathrm{E}_{ 0})
\biggl(\exp\biggl[\frac{\gamma}{\kappa} \int_0^t ds\,
\xi\biggl(X(s),\frac{t-s}\kappa\biggr) \biggr] \biggr),
\end{equation}
where $X=X^1$.
For large $\kappa$, the $\xi$-field in (\ref{lambda*scal}) evolves
slowly and
therefore does not manage to cooperate with the $X$-process in
determining the
growth rate. As a result, the expectation over the $\xi$-field can be computed
via a \textit{Gaussian approximation}, which we expect to become sharp
in the limit
as $\kappa\to\infty$, that is,
%
\begin{eqnarray}\label{Gappr1}
\quad &&\Lambda_1^*(\kappa;t) - \frac{\rho\gamma}{\kappa}\nonumber\\
\quad &&\qquad = \frac{1}{t} \log(\mathbb{E}_{ \mu_\rho}\otimes\mathrm{E}_{ 0})
\biggl(\exp\biggl[\frac{\gamma}{\kappa} \int_0^t \,ds
\biggl(\xi\biggl(X(s),\frac{t-s}\kappa\biggr)-\rho\biggr) \biggr] \biggr)\nonumber\\
\quad &&\qquad  \approx\frac{1}{t} \log
\mathrm{E}_{ 0} \biggl(\exp\biggl[\frac{\gamma^2}{2\kappa^2}
\int_0^t ds \int_0^t du\,
\mathbb{E}_{ \mu_\rho} \biggl(
\biggl[\xi\biggl(X(s),\frac{t-s}\kappa\biggr)-\rho\biggr]\\
\quad &&\qquad \quad \hphantom{\frac{1}{t} \log
\mathrm{E}_{ 0} \biggl(\exp\biggl[\frac{\gamma^2}{2\kappa^2}
\int_0^t ds \int_0^t du\,
\mathbb{E}_{ \mu_\rho} \biggl(
\biggl[}{}\times
\biggl[\xi\biggl(X(u),\frac{t-u}\kappa\biggr)-\rho\biggr] \biggr) \biggr]
\biggr).\nonumber\hspace*{-12pt}
\end{eqnarray}
(In essence, what happens here is that the asymptotics for $\kappa\to
\infty$ is
driven by moderate deviations of the $\xi$-field, which fall in the
Gaussian regime.)
Next, by Lemma \ref{2corr}, for any $0\leq s\leq u\leq t$ we have
%
\begin{eqnarray}\label{Gappr3}
&&\mathbb{E}_{ \mu_\rho} \biggl( \biggl[\xi\biggl(X(s),\frac{t-s}\kappa\biggr)
-\rho\biggr] \biggl[\xi\biggl(X(u),\frac{t-u}\kappa\biggr)-\rho\biggr] \biggr)\nonumber\\[-8pt]\\[-8pt]
&&\qquad  = C \int_0^\infty dv\, p_{(u-s)/\kappa+v} (X(s),X(u) ),\nonumber
\end{eqnarray}
where $C=\rho(1-\rho)/G_d$. Hence,
%
\begin{eqnarray}\label{limfinal}
\lim_{\kappa\to\infty} 2d\kappa[\lambda_1(\kappa)-\rho\gamma]
&=& \lim_{\kappa\to\infty}
2d\kappa^2 \biggl[\lambda_1^*(\kappa)-\frac{\rho\gamma}{\kappa} \biggr]\nonumber\\
&=& \lim_{\kappa\to\infty}2d\kappa^2\lim_{t\to\infty}
\biggl[\Lambda_1^*(\kappa;t)-\frac{\rho\gamma}{\kappa} \biggr]\\
&=& \lim_{\kappa\to\infty}
2d\kappa^2\lim_{t\to\infty} I(\kappa;t)\nonumber
\end{eqnarray}
with
%
\begin{eqnarray}\label{Gappr7}
\qquad &&I(\kappa;t)\nonumber\\
\qquad &&\qquad =\frac{1}{t}\log\mathrm{E}_{ 0} \biggl(
\exp\biggl[\frac{C\gamma^2}{\kappa^2}\int_0^t ds \int_s^t du
\int_0^\infty dv \,p_{(u-s)/{\kappa}+v} (X(s),X(u) ) \biggr]
\biggr)\\
\qquad &&\qquad \approx\frac{C\gamma^2}{t\kappa^2}\int_0^t ds \int_s^t du
\int_0^\infty dv\,
\mathrm{E}_{ 0} \bigl(p_{(u-s)/{\kappa}+v} (X(s),X(u) ) \bigr).\nonumber
\end{eqnarray}
In the last line of (\ref{Gappr7}), a \textit{linear approximation} is
made in the
expectation over the random walk $X$, which we expect to become sharp
in the limit
as $\kappa\to\infty$ in $d \geq6$. Next, for any $0\leq s\leq u\leq
t$ and $T\geq0$,
%
\begin{eqnarray}\label{Gappr10}
\mathrm{E}_{ 0} (p_T(X(s),X(u)) )
&=&\sum_{x,y\in\mathbb{Z}^d} p_{2ds}(0,x) p_{2d(u-s)}(x,y) p_T(x,y)\nonumber\\
&=&\sum_{x\in\mathbb{Z}^d} p_{2ds}(0,x) p_{2d(u-s)+T}(x,x)\\
&=&p_{2d(u-s)+T}(0,0).\nonumber
\end{eqnarray}
Here, we use that $p(\cdot,\cdot)$ is a simple random walk, so that
$\xi$ fits with
$X$. We therefore have
%
\begin{equation}\label{Gappr16}
\mbox{r.h.s. (\ref{Gappr7})}
=\frac{C\gamma^2}{t\kappa^2} \int_0^t ds \int_s^t du \int
_0^\infty dv\,
p_{2d(u-s)1[\kappa]+v} (0,0 ),
\end{equation}
where we abbreviate $1[\kappa] = 1 + \frac{1}{2d\kappa}$. Rewriting
%
\begin{eqnarray}\label{Gappr24}
&&\frac1t\int_0^t ds \int_s^t du
\int_0^\infty dv\,
p_{2d(u-s)1[\kappa]+v} (0,0 )\nonumber\\[-8pt]\\[-8pt]
&&\qquad  =\int_0^t dw \int_0^\infty dv
\biggl(\frac{t-w}{t} \biggr)
p_{2dw1[\kappa]+v} (0,0 ),\nonumber
\end{eqnarray}
we get from (\ref{Gappr7})--(\ref{Gappr16}) that
%
\begin{eqnarray}\label{Gappr28}
\lim_{t\rightarrow \infty}I(\kappa;t)
&=&\frac{C\gamma^2}{2d\kappa^21[\kappa]}
\int_0^\infty dw \int_0^\infty dv\,
p_{w+v}(0,0)\nonumber\\[-8pt]\\[-8pt]
&=&\frac{C\gamma^2}{2d\kappa^21[\kappa]} G_d^\ast.\nonumber
\end{eqnarray}
Recalling (\ref{limfinal}), we arrive at $(\ref{limlamb*})$ for $d
\geq6$.

3. \textit{Polaron term}.
Where does the term with $\mathcal{P}_5$ come from? We expect this
term to
arise from
the part of the integral in the exponent in the first line of (\ref{Gappr7})
with $(u-s)/\kappa$ and $v$ of order $\kappa^2$, as we will argue
next. Put
$\mathbb{Z}^d_\kappa=\kappa^{-1}\mathbb{Z}^d$ and, for $t \geq0$
and $x,y\in\mathbb{Z}
^d_\kappa$, define
%
\begin{equation}\label{kscaldefs}
X^\kappa(t) = \kappa^{-1} X(\kappa^2t),\qquad
p_t^\kappa(x,y) = \kappa^d p_{2d\kappa^2 t}(\kappa x,\kappa y).
\end{equation}
In the limit as $\kappa\to\infty$, $(X^\kappa(t))_{t\geq0}$
converges weakly to
Brownian motion, while $(p_t^\kappa(\cdot,\cdot))_{t\geq0}$
converges to the
corresponding family of Gaussian transition kernels $(p_t^G(\cdot
,\cdot))_{t\geq0}$
given by
%
\begin{equation}\label{pGdef}
p_t^G(x,y) = (4\pi t)^{-d/2} \exp[-\|x-y\|^2/4t],\qquad  x,y\in\mathbb{R}^d.
\end{equation}
After scaling, the part we are after is approximately
%
\begin{equation}\label{intscal1}
\qquad\hspace*{8pt} C\gamma^2 \kappa^{4-d} \int_0^{\kappa^{-2}t} ds
\int_{s+\varepsilon\kappa}^{s+K\kappa} du
\int_0^K dv\, p^G_{1/(2d) ((u-s)/{\kappa}+v )} (X^\kappa
(s),X^\kappa(u) ),\hspace*{-8pt}
\end{equation}
where $0 < \varepsilon\ll1 \ll K <\infty$. For $\delta>0$, divide
the first and the second integral
in (\ref{intscal1}) into pieces of length $\delta\kappa$, and define
the occupation
time measures
%
\begin{equation}\label{occmeasdef}
\Xi_w^\kappa(A) = \frac{1}{\delta\kappa} \int_w^{w+\delta\kappa
} 1_A(X^\kappa(u)) \,du,\qquad
w \geq0, A \subset\mathbb{R}^d \mbox{ Borel}.
\end{equation}
Then, when $\delta\ll\varepsilon$, $(u-s)/\kappa$ is almost
constant on time intervals
of length $\delta\kappa$ and, consequently,
%
\begin{eqnarray}\label{intscal2}
\qquad  \mbox{(\ref{intscal1})} &\approx& C\gamma^2 \kappa^{4-d}
\int_0^{\kappa^{-2}t} ds \int_{s+\varepsilon\kappa}^{s+K\kappa} du
\int_0^K dv\nonumber\\[-8pt]\\[-8pt]
&&\hphantom{C\gamma^2 \kappa^{4-d}
\int_0^{\kappa^{-2}t}}
{}\times \int_{\mathbb{R}^d} \Xi_s^\kappa(dx) \int_{\mathbb
{R}^d} \Xi
_u^\kappa(dy) p_{1/{(2d)} ((u-s)/{\kappa}+v )}^G(x,y).\nonumber
\end{eqnarray}
Using the large deviation principle for $\Xi_{(\cdot)}^\kappa$ as
$\kappa\to\infty$,
we find that the contribution of (\ref{intscal2}) to $I(\kappa;t)$
for large
$\kappa$ is approximately
%
\begin{eqnarray}\label{varform1}
\qquad &&\frac{1}{t} \sup_{\mu_{(\cdot)}} \biggl[C\gamma^2\kappa^{4-d} \int
_0^{\kappa^{-2}t} ds
\int_{s+\varepsilon\kappa}^{s+K\kappa} du \int_0^K dv \int
_{\mathbb{R}^d}
\mu_s(dx) \int_{\mathbb{R}^d} \mu_u(dy) \nonumber\\[-8pt]\\[-8pt]
\qquad &&\hspace*{102pt}
{}\times p_{1/({2d}) ((u-s)/{\kappa}+v )}^G(x,y)
- \int_0^{\kappa^{-2}t} J(\mu_s) \,ds \biggr],\nonumber
\end{eqnarray}
where the supremum is taken over all probability measure-valued paths
$\mu_{(\cdot)}$ and
%
\begin{equation}\label{Iiden}
J(\nu) = \cases{
\bigl\|\nabla\sqrt{d\nu/d\lambda} \bigr\|_2^2, &\quad if $
\nu\ll\lambda$,\cr
\infty,&\quad  otherwise,
}
\end{equation}
with $\lambda$ the Lebesgue measure on $\mathbb{R}^d$. By the
convexity of the
large deviation rate function $J$, the supremum in (\ref{varform1})
diagonalizes
and reduces to
%
\begin{eqnarray}\label{varform2}
\mbox{(\ref{varform1})} &=& \frac{1}{\kappa^2} \sup_{\nu}
\biggl[C \gamma^2\kappa^{4-d} \int_{\mathbb{R}^d} \nu(dx) \int_{\mathbb
{R}^d} \nu(dy)
\int_{\varepsilon\kappa}^{K\kappa} du \int_0^K
dv\nonumber\\[-8pt]\\[-8pt]
&&\hspace*{106pt}
{}\times p_{1/({2d})(u/{\kappa}+v)}^G(x,y) - J(\nu) \biggr].\nonumber
\end{eqnarray}
Putting $u=2d\kappa\tilde{u}$, $v=2d\tilde{v}$ and letting
$\varepsilon\downarrow0$
and $K\to\infty$, we end up with a contribution to
$\lim_{\kappa\to\infty}2d\kappa^2$ $\lim_{t\to\infty} I(\kappa;t)$
of the form
%
\begin{eqnarray}\label{varform3}
&&2d \sup_{\nu} \biggl[(2d)^2 C\gamma^2\int_{\mathbb{R}^d} \nu(dx) \int
_{\mathbb{R}^d}
\nu(dy)
\int_0^\infty d\tilde{u} \int_0^\infty d\tilde{v}\, p_{\tilde
{u}+\tilde{v}}^G(x,y)\nonumber\\[-8pt]\\[-8pt]
&&\hspace*{255pt}
{} - J(\nu) \biggr]\nonumber
\end{eqnarray}
in $d=5$ and zero in $d \geq6$. In $d=5$ we have from (\ref{pGdef})
%
\begin{equation}\label{Gint}
\int_0^\infty d\tilde{u} \int_0^\infty d\tilde{v}\, p_{\tilde
{u}+\tilde{v}}^G(x,y)
= \int_0^\infty dt\, t p_t^G(x,y) = \frac{1}{16\pi^2\|x-y\|}.
\end{equation}
Substituting this into (\ref{varform3}), putting $\nu=f^2\lambda$
and recalling (\ref{Iiden}),
we get
%
\begin{equation}\label{varform4}
\hspace*{18pt}\mbox{(\ref{varform3})} =
2d \sup_{\|f\|_2=1} \biggl[(2d)^2 C\gamma^2\int_{\mathbb{R}^5}\int
_{\mathbb{R}^5} dx \,dy\,
\frac{f^2(x) f^2(y)}{16\pi^2\|x-y\|} - \|\nabla f\|_2^2 \biggr].\hspace*{-18pt}
\end{equation}
Scaling of $f$ shows that the supremum with the prefactor
$(2d)^2C\gamma^2$ equals
$((2d)^2C\gamma^2)^2$ times the supremum without this prefactor.
Hence, we get
%
\begin{equation}\label{varform5}
(\ref{varform3}) = 2d ((2d)^2C\gamma^2 )^2 \mathcal{P}_5,
\end{equation}
where we recall (\ref{P5def}). This is precisely the ``polaron-type''
term in (\ref{limlamb*})
for $p=1$.

The heuristic argument in parts 2 and 3 follows a line of thought that
was made rigorous in
G\"artner and den Hollander \cite{garhol06} and G\"artner, den
Hollander and Maillard
\cite{garholmai07,garholmai08pr} for the case where $\xi$ is
a field of independent
simple random walks in a Poisson equilibrium, respectively, a simple
symmetric exclusion
process in a Bernoulli equilibrium. We refer to these papers for
further details. There it
is also explained why for $p\in\mathbb{N}\setminus\{1\}$ the
polaron term
is $p^2$ times that for $p=1$.
\end{appendix}

%

\printaddresses

\end{document}